\documentclass[12pt,a4paper]{article}

\usepackage[round,semicolon]{natbib} 
\setcitestyle{notesep={, },yysep={,},aysep={,},citesep={,}}

\usepackage{anysize}
\usepackage{color}
\usepackage{graphicx}
\usepackage[active]{srcltx} 

\usepackage{amsmath,amssymb,amsthm,amsopn}
\usepackage{verbatim,graphicx}

\RequirePackage[latin1]{inputenc}
\RequirePackage[english]{babel}
\RequirePackage{amssymb} \RequirePackage{amsfonts}
\RequirePackage{amsmath} \RequirePackage{amsthm}
\RequirePackage{latexsym}

\newtheorem{theorem}{Theorem}

\newtheorem{remark}{Remark}

\newcommand{\R}{\mathbb{R}}
\newcommand{\N}{\mathbb{N}}

\def\MISE{\mathrm{MISE}}

\def\E{\mathrm{E}}
\def\Var{\mathrm{Var}}

\def\tF{\tilde{F}}
\def\bK{\bar{K}}
\def\bF{\bar{F}}

\newcommand{\halmoseq}{\tag*{$\blacksquare$}}

\theoremstyle{plain}

\title{A note on boundary kernels for distribution function estimation}

\author{
Carlos Tenreiro\footnote{CMUC, Department of Mathematics, University of Coimbra, Apartado 3008, 3001--501 Coimbra, Portugal.
E-mail: tenreiro@mat.uc.pt. URL: \texttt{http://www.mat.uc.pt/$\sim$tenreiro/}}
}

\date{January 17, 2015}

\begin{document}


\maketitle

\begin{abstract}
\noindent
The use of second order boundary kernels for distribution function estimation was recently addressed in the literature
(C.~Tenreiro, 2013, Boundary kernels for distribution function estimation, \textit{REVSTAT--Statistical Journal}, 11, 169--190). In this note we return to the subject by considering an enlarged class of boundary kernels that shows it self to be especially performing when the classical kernel distribution function estimator suffers from severe boundary problems.

\bigskip
\noindent {\sc Keywords:} Distribution function estimation; kernel estimator; boundary kernels.

\bigskip
\noindent AMS 2010 {\sc subject classifications:} 62G05, 62G20

\end{abstract}

\newpage

\section{Introduction}

Given $X_1,\dots,X_n$ independent copies of an absolutely continuous real random variable
with unknown density and distribution functions $f$ and $F$, respectively, the classical kernel estimator of
$F$ introduced by authors such as
\citet{Tia:63}, \citet{Nad:64} or \citet{WatL:64}, is defined, for $x\in\R$, by
\begin{equation} \label{classical}
\bF_{nh}(x)=\frac{1}{n} \sum_{i=1}^{n} \bK \left( \frac{x-X_i}{h} \right),
\end{equation}
where, for $u\in\R$,
$$\bK(u) = \int_{-\infty}^u K(v)dv,$$
with $K$ a kernel on $\R$, that is, a bounded and symmetric probability density function with support $[-1,1]$
and $h=h_n$ a sequence of strictly
positive real numbers converging to zero when $n$ goes to infinity.
For some recent references on this classical estimator see
\citet{GinN:09}, \citet{ChaR:10}, \citet{MasS:11} and \citet{ChaMT:14}.

If the support of $f$ is known to be the finite interval $[a,b]$,
the previous kernel estimator suffers from boundary problems
if $F_+^\prime(a) \neq 0$ or $F_-^\prime(b) \neq 0$.
This question is addressed in \citet{Ten:13} by extending to the distribution function estimation framework
the approach followed in nonparametric regression and density function estimation
by authors such as \citet{GasM:79}, \citet{Ric:84}, \citet{GasMM:85} and \citet{Mul:91}.
Specially, the author considers the boundary modified kernel distribution function estimator given by
\begin{equation} \label{boundary}
\tF_{nh}(x) =
\left\{
\begin{array}{ll}
0, & x \leq a \\
{\displaystyle \frac{1}{n} \sum_{i=1}^{n} \bK_{x,h} \left( \frac{x-X_i}{h} \right)}, & a < x < b \\
1, & x \geq b, \\
\end{array}
\right.
\end{equation}
where $0<h\leq (b-a)/2$ and
$$\bK_{x,h}(u)=
\left\{
\begin{array}{ll}
\bK^L(u;(x-a)/h), & a < x < a+h \\
\bK(u), & a+h \leq x \leq b-h \\
\bK^R(u;(b-x)/h), & b-h < x < b, \\
\end{array}
\right.
$$
with
$$\bK^L(u;\alpha)=\int_{-\infty}^u K^L(v;\alpha) dv \quad \mbox{and} \quad
\bK^R(u;\alpha)=1-\int_u^{+\infty} K^R(v;\alpha) dv,$$
where $K^L(\cdot;\alpha)$ and $K^R(\cdot;\alpha)$
are, respectively, left and right boundary kernels for $\alpha\in \,]0,1[$, that is,
their supports are contained in the intervals $[-1,\alpha]$ and $[-\alpha,1]$, respectively,
and $|\mu_{0,\ell}|(\alpha)=\int |K^\ell(u;\alpha)| \, du<\infty$ for all $\alpha\in\,]0,1[$ and $\ell=L,R$
(here and bellow integrals without integrations limits are meant over the whole real line).

For ease of presentation, from now on we assume that the right boundary kernel $K^R$
is given by $K^R(u;\alpha)=K^L(-u;\alpha)$, the reason why only the left boundary kernel is
mentioned in the following discussion.
By assuming that $K^L(\cdot;\alpha)$ is a second order kernel, that is,
\begin{equation} \label{C1}
\mu_{0,L}(\alpha) = 1, \; \mu_{1,L}(\alpha) = 0 \mbox{ and }
\mu_{2,L}(\alpha) \neq 0, \mbox{ for all }\alpha\in\,]0,1[,
\end{equation}
where we denote
$$\mu_{k,L}(\alpha)=\int u^k K^L(u;\alpha) \, du, \mbox{ for } k\in\N,$$
\citet{Ten:13} shows that the previous estimator
is free of boundary problems and
that the theoretical advantage of using boundary kernels is compatible
with the natural property of getting a proper distribution
function estimate. In fact, it is easy to see that the kernel distribution function estimator based on each one
of the second order left boundary kernels
\begin{equation} \label{K1}
K^L_1(u;\alpha) = (2\bK(\alpha)-1)^{-1} K(u)I(-\alpha \leq u \leq \alpha),
\end{equation}
where we assume that $K$ is such that $\int_0^\alpha K(u) du >0$ for all $\alpha>0$,
and
\begin{equation} \label{K2}
K^L_2(u;\alpha) = K(u/\alpha)/\alpha,
\end{equation}
is, with probability one, a continuous probability distribution function \citep[see][Examples 2.2 and 2.3]{Ten:13}.
Additionally, the author shows that the Chung-Smirnov law of iterated logarithm
is valid for the new estimator and has presented
an asymptotic expansion for its mean integrated squared error, from which the
choice of $h$ is discussed \citep[see][Theorems 3.2, 4.1 and 4.2]{Ten:13}.

\begin{figure}[!t]

\centering
\begin{tabular}{c@{\hspace{2em}}c}
\multicolumn{2}{c}{\small \textsf{\hspace{2em}$K^L_1(u;\alpha)$}}  \\*[-0.0ex]
\includegraphics[scale=0.4]{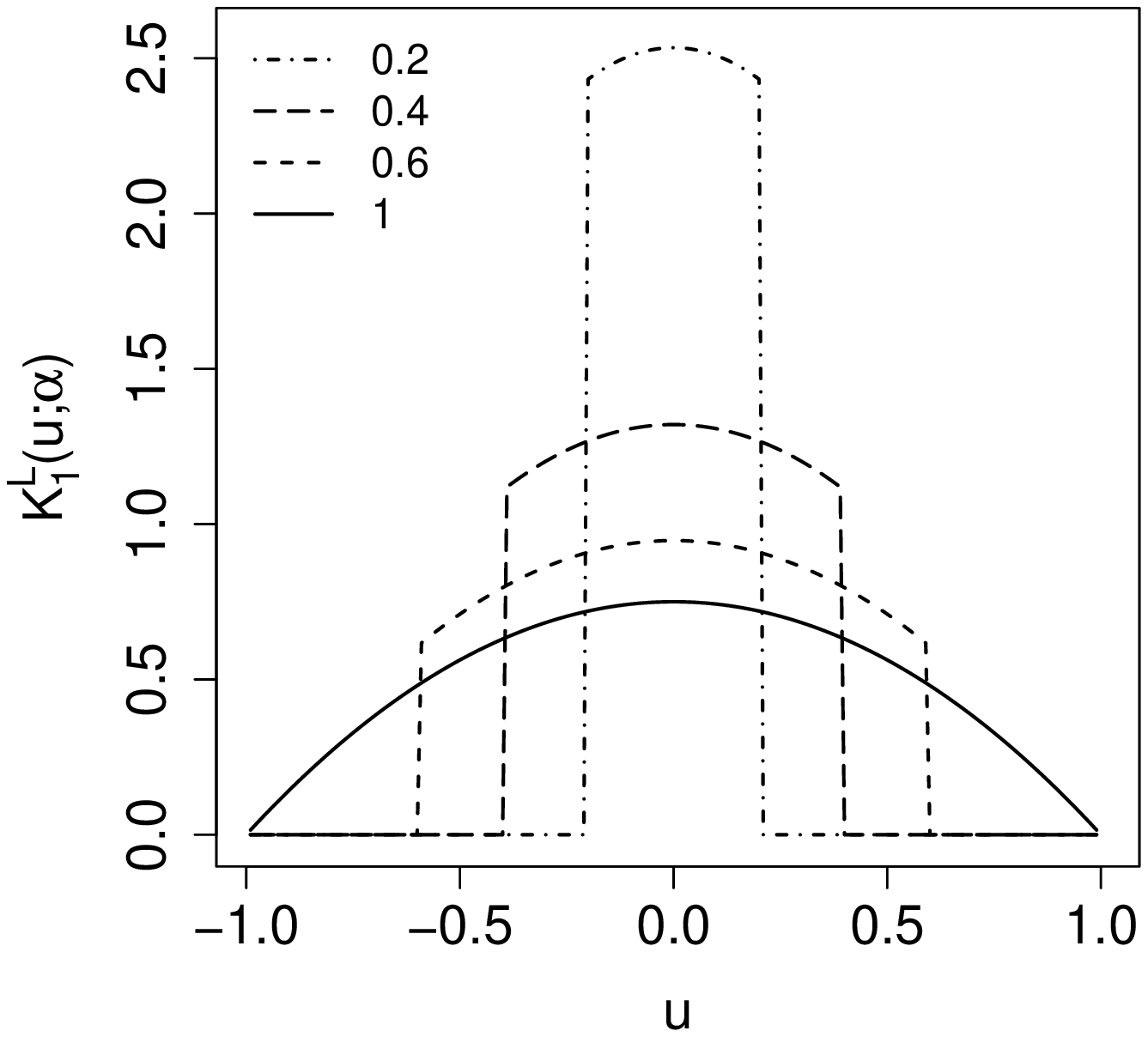} & \includegraphics[scale=0.4]{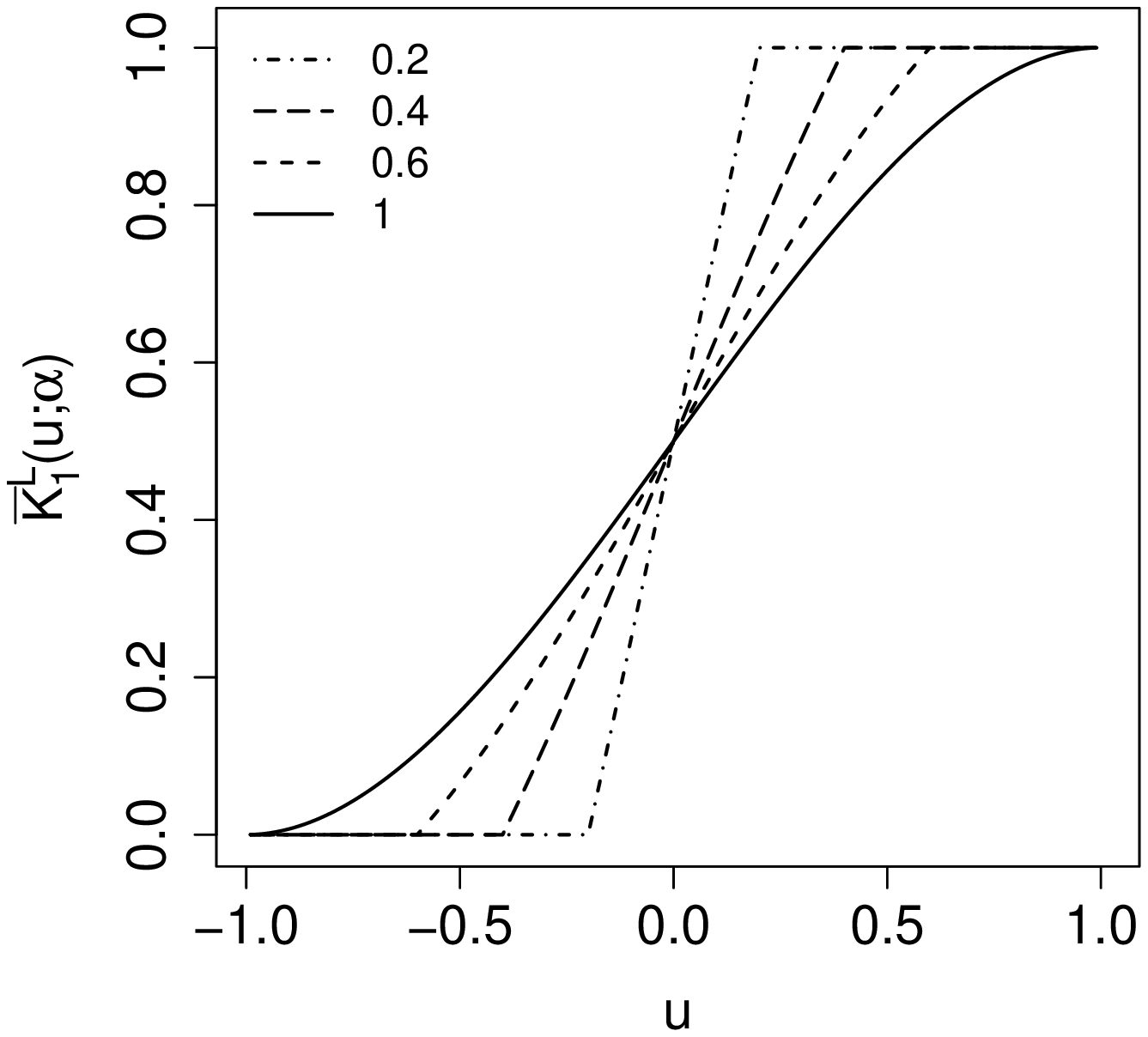} \\*[-1ex]
\multicolumn{2}{c}{\small \textsf{\hspace{2em}$K^L_2(u;\alpha)$}}  \\
\includegraphics[scale=0.4]{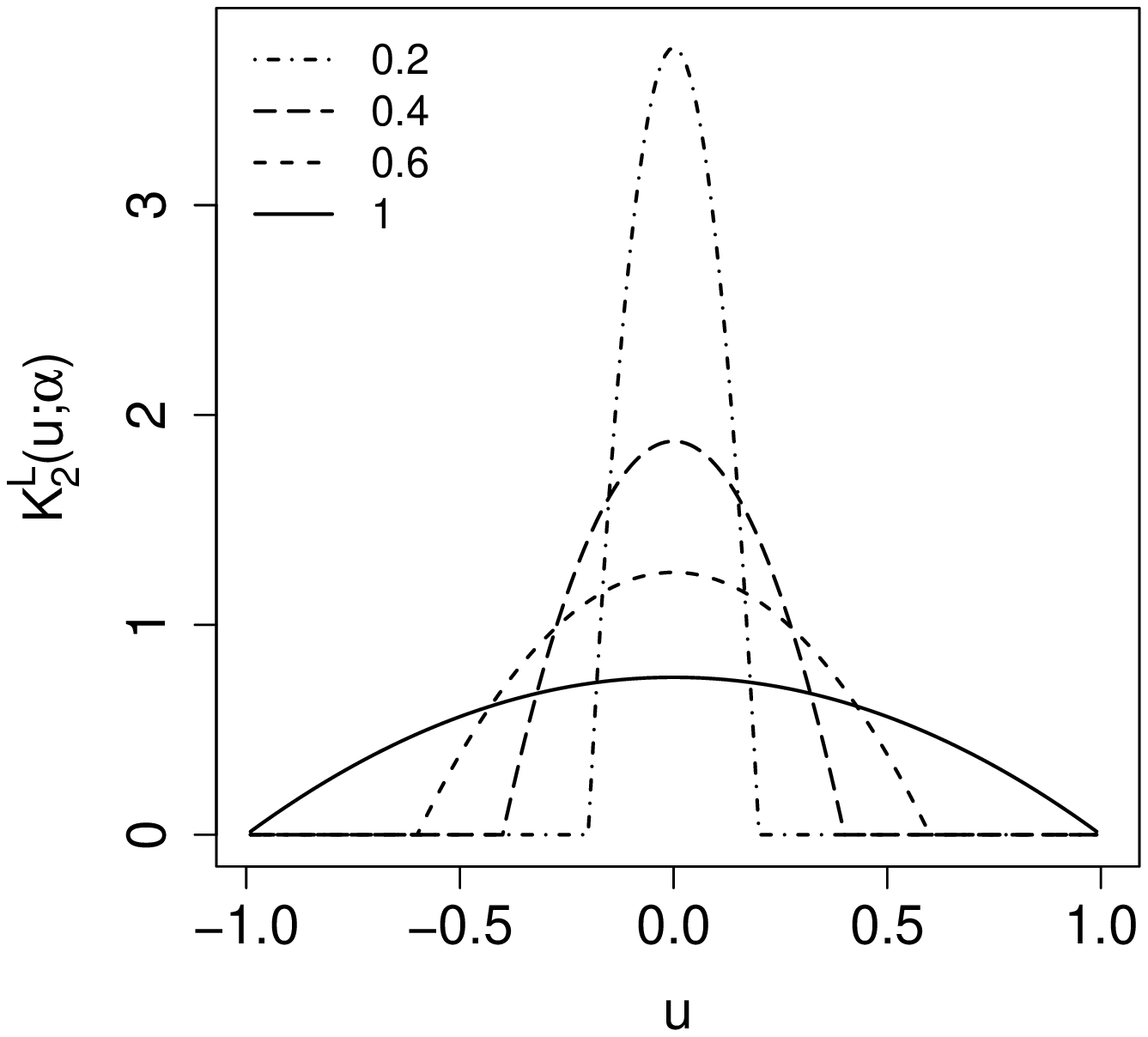} & \includegraphics[scale=0.4]{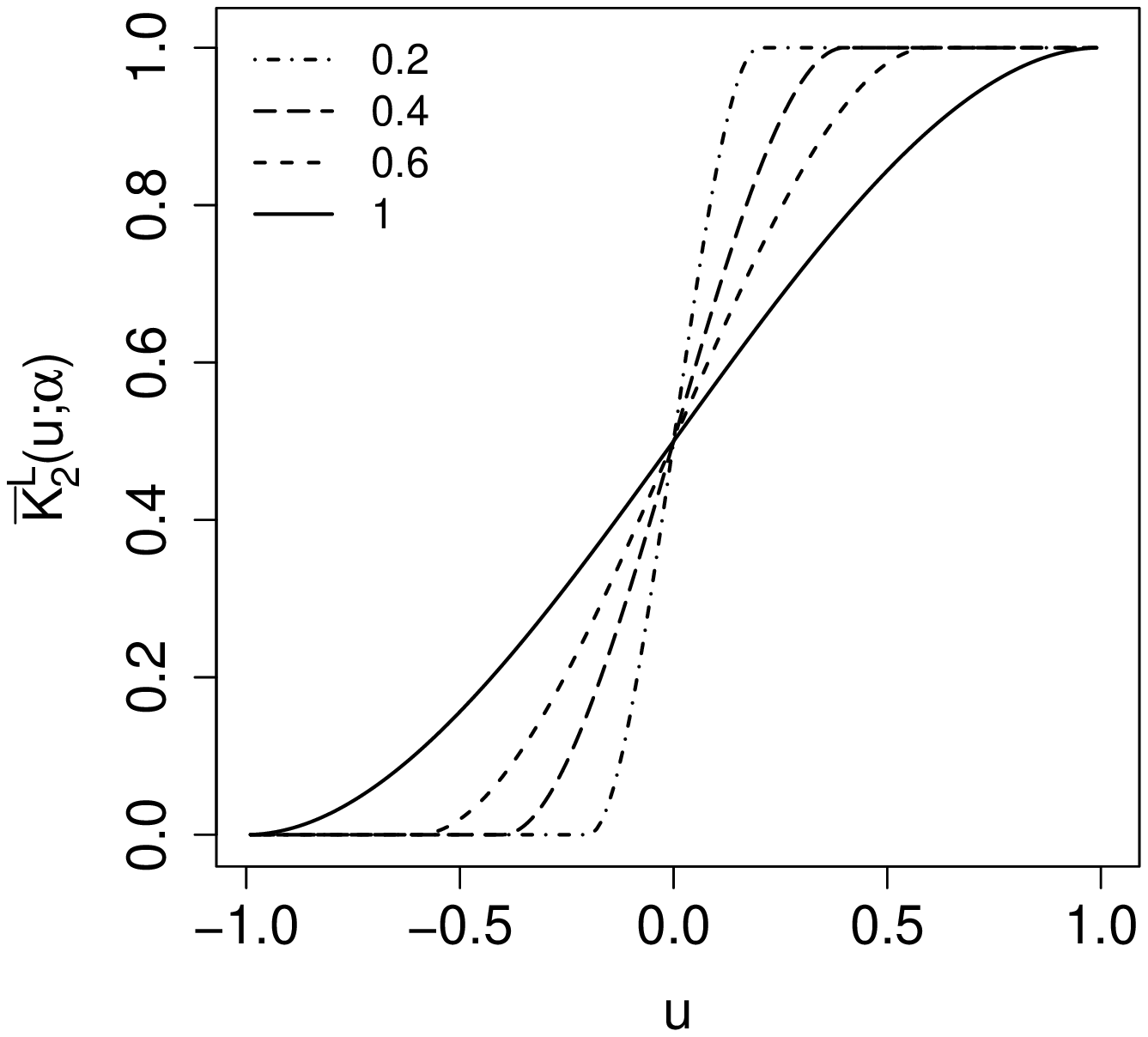} \\*[-1ex]
\multicolumn{2}{c}{\small \textsf{\hspace{2em}$K^L_3(u;\alpha)$}}  \\
\includegraphics[scale=0.4]{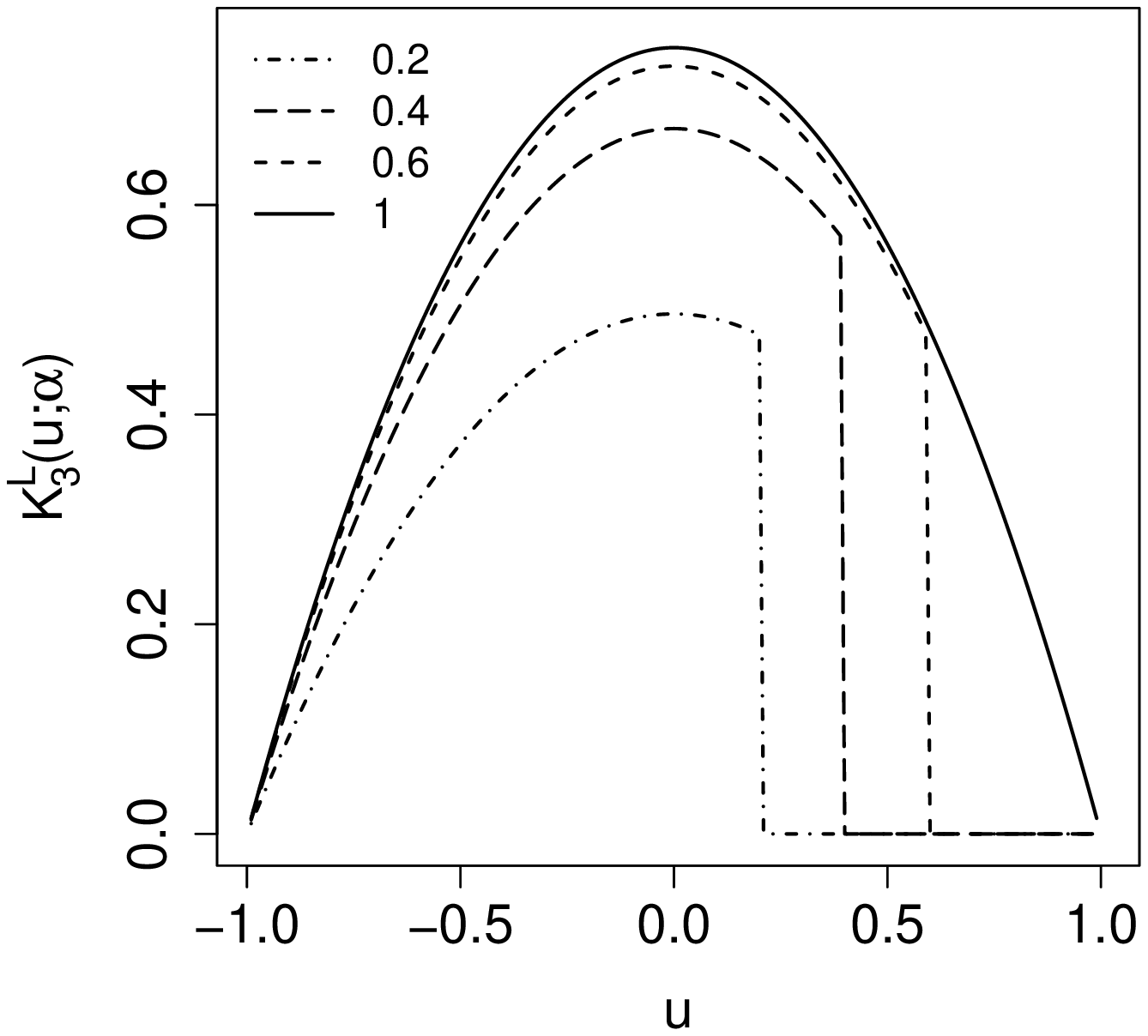} & \includegraphics[scale=0.4]{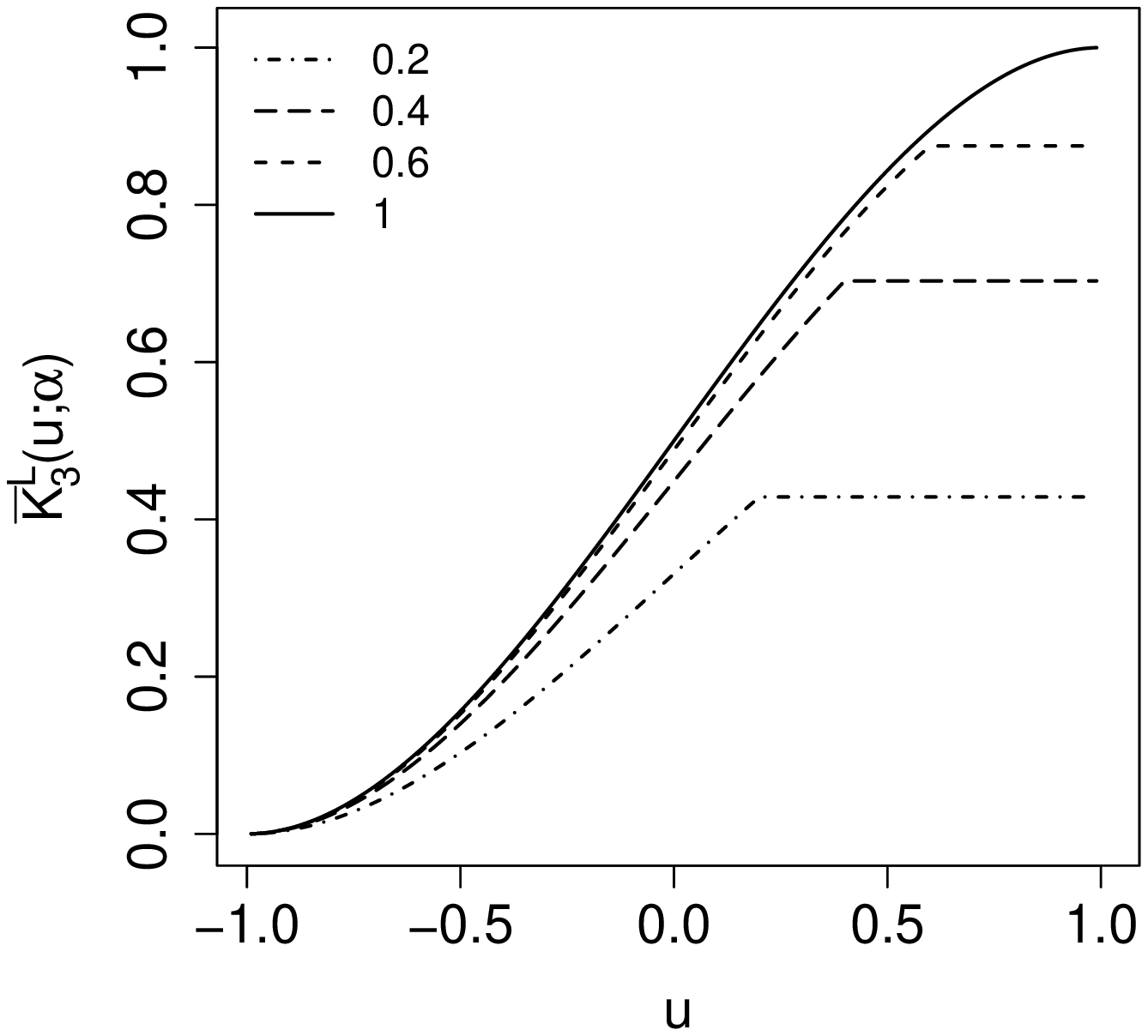}
\end{tabular}

\vspace{-1ex}
\caption{\it Left boundary kernels $K^L_q(u;\alpha)$ (left column) and $\bK^L_q(u;\alpha)$ (right column) for $q=1,2,3$, where $K$ is the Epanechnikov kernel $K(t)=\frac{3}{4}(1-t^2) I(|t|\leq 1)$.} \label{kernels}
\end{figure}

A careful analysis of the asymptotic expansions presented in
\citet[][p.~171, 178]{Ten:13} for the local bias and
the integrated squared bias of estimator (\ref{classical}), suggests that the previous properties may still be valid for
all the boundary kernels satisfying the less restricted condition
\begin{equation} \label{C2}
\alpha \left( 1 - \mu_{0,L}(\alpha) \right) + \mu_{1,L}(\alpha) = 0,\mbox{ for all }\alpha\in\,]0,1[,
\end{equation}
which is in particular fulfilled by the left boundary kernel
\begin{equation} \label{K3}
K^L_3(u;\alpha) = \alpha K(u) I( -1\leq u \leq \alpha) \big/ (\alpha \mu_{0,\alpha}(K) - \mu_{1,\alpha}(K)),
\end{equation}
where we denote
$\mu_{k,\alpha}(K)=\int_{-1}^\alpha u^k K(u) \, du$, for $k\in\N$ (see Figure \ref{kernels}).
If $K$ is a continuous density function, it is not hard to prove that the kernel distribution function estimator
based on this left boundary kernel is, with probability one, a continuous probability distribution function.

The main purpose of this note is to show that the results presented in \citet{Ten:13} for the class of second order
boundary kernels are still valid for the enlarged class of boundary kernels that satisfy assumption (\ref{C2}).
This objective is achieved in Sections \ref{secboundary} and \ref{global} where we study the boundary
and global behaviour of the boundary modified kernel distribution function estimator $\tF_{nh}$.
In Section \ref{compar} we present exact finite sample comparisons between the distribution
function kernel estimators based on the left boundary kernels $K^L_q(u;\alpha)$, for $q=1,2,3$,
given by (\ref{K1}), (\ref{K2}) and (\ref{K3}), respectively.
We conclude that the boundary kernel $K_3^L$ is especially performing when the
classical kernel estimator suffers from severe boundary problems.
All the proofs can be found in Section \ref{proofs}.
The plots and simulations in this paper were carried out using the \textsf{R} software \citep{R:11}.

\section{Boundary behaviour} \label{secboundary}

In this section we study the boundary behaviour of the kernel distribution function estimator $\tF_{nh}(x)$ by presenting
asymptotic expansions for its bias and variance with $x$ in the boundary region.
We will restrict our attention to the left boundary region $]a,a+h[$. However, similar similar results
are valid for the right boundary region $]b-h,b[$.

\begin{theorem} \label{theoremboundL}
If $K^L(u;\alpha)$ satisfies condition (\ref{C2}) with
$$\sup_{\alpha\in\, ]0,1[} |\mu_{0,L}|(\alpha) < \infty,$$
and the restriction of $F$ to the interval $[a,b]$ is twice continuously differentiable, we have:

\medskip

a) $$\sup_{x\in\,]a,a+h[} \left| \E \tF_{nh}(x) - F(x) - \frac{h^2}{2} F^{\prime\prime}(x) \mu_{L}\big((x-a)/h\big) \right| = o(h^2).$$
where
$$\mu_L(\alpha)=\mu_{2,L}(\alpha) - \alpha \mu_{1,L}(\alpha), \; \alpha\in\,]0,1[;$$

b) $$\sup_{x\in\,]a,a+h[} \left| \Var \tF_{nh}(x) - \frac{F(x)\big(1-F(x)\big)}{n}
+ \frac{h}{n}\, F^\prime(x) \nu_L\big((x-a)/h\big) \right| = O(n^{-1}h^2),$$
where
$$\nu_L(\alpha)=m_{1,L}(\alpha)+\alpha(1-\mu_{0,L}(\alpha)^2), \; \alpha\in\,]0,1[,$$
with
$m_{1,L}(\alpha)=\int  uB^L(u;\alpha)\,du,$
and $B^L(u;\alpha)=2\bK^L(u;\alpha)K^L(u;\alpha)$.
\end{theorem}

\medskip

\begin{remark}
{\rm
The previous expansions for the bias and variance of $\tF_{nh}(x)$ extend those
presented in \citet[p.~174]{Ten:13} for second order boundary kernels, in which case
$\mu_L(\alpha)=\mu_{2,L}(\alpha)$ and $\nu_L(\alpha)=m_{1,L}(\alpha)$, for $\alpha\in\,]0,1[$.
}
\end{remark}

\begin{figure}[!t]

\centering
\begin{tabular}{c@{\hspace{2em}}c}
\includegraphics[scale=0.45]{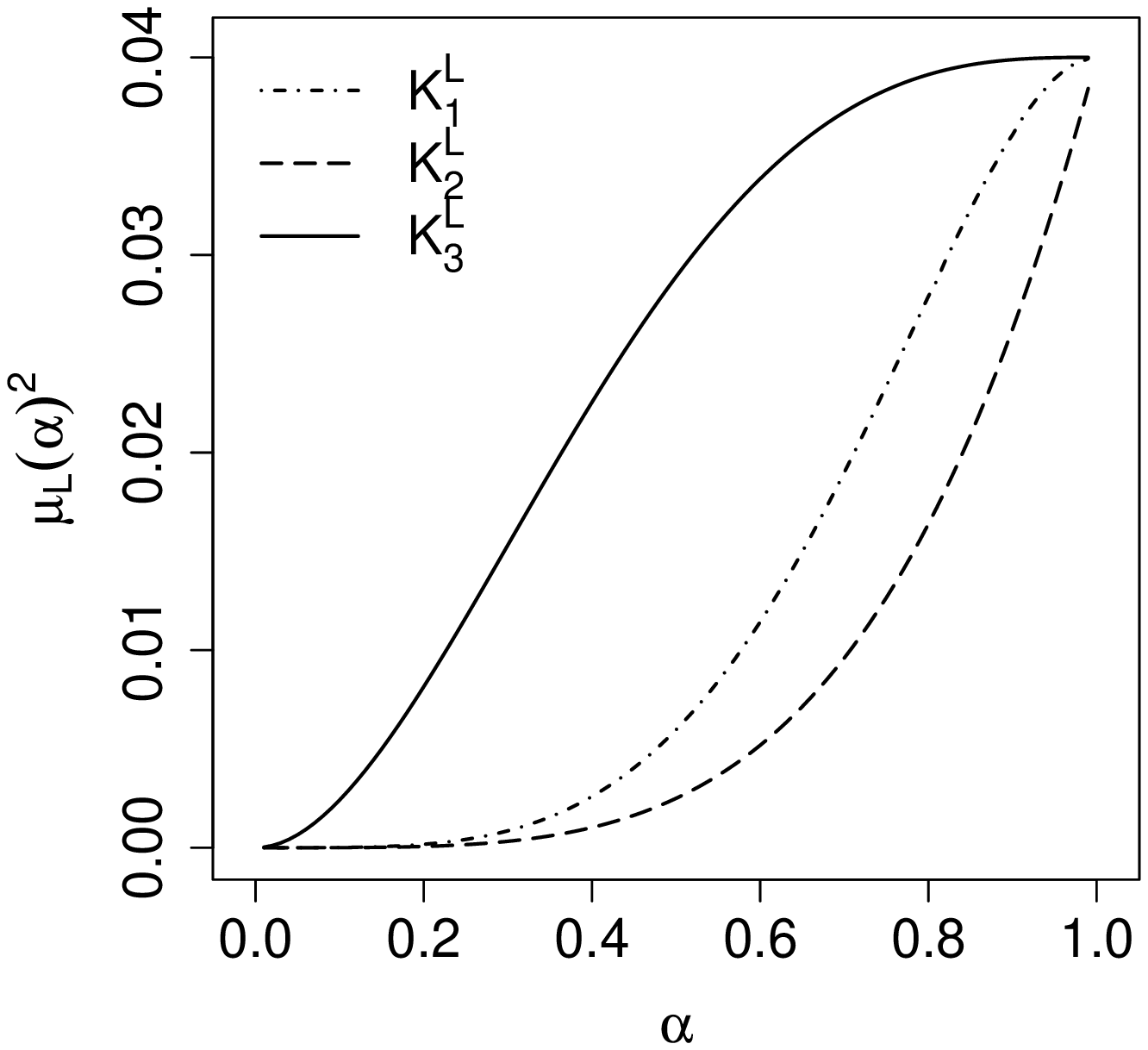} & \includegraphics[scale=0.45]{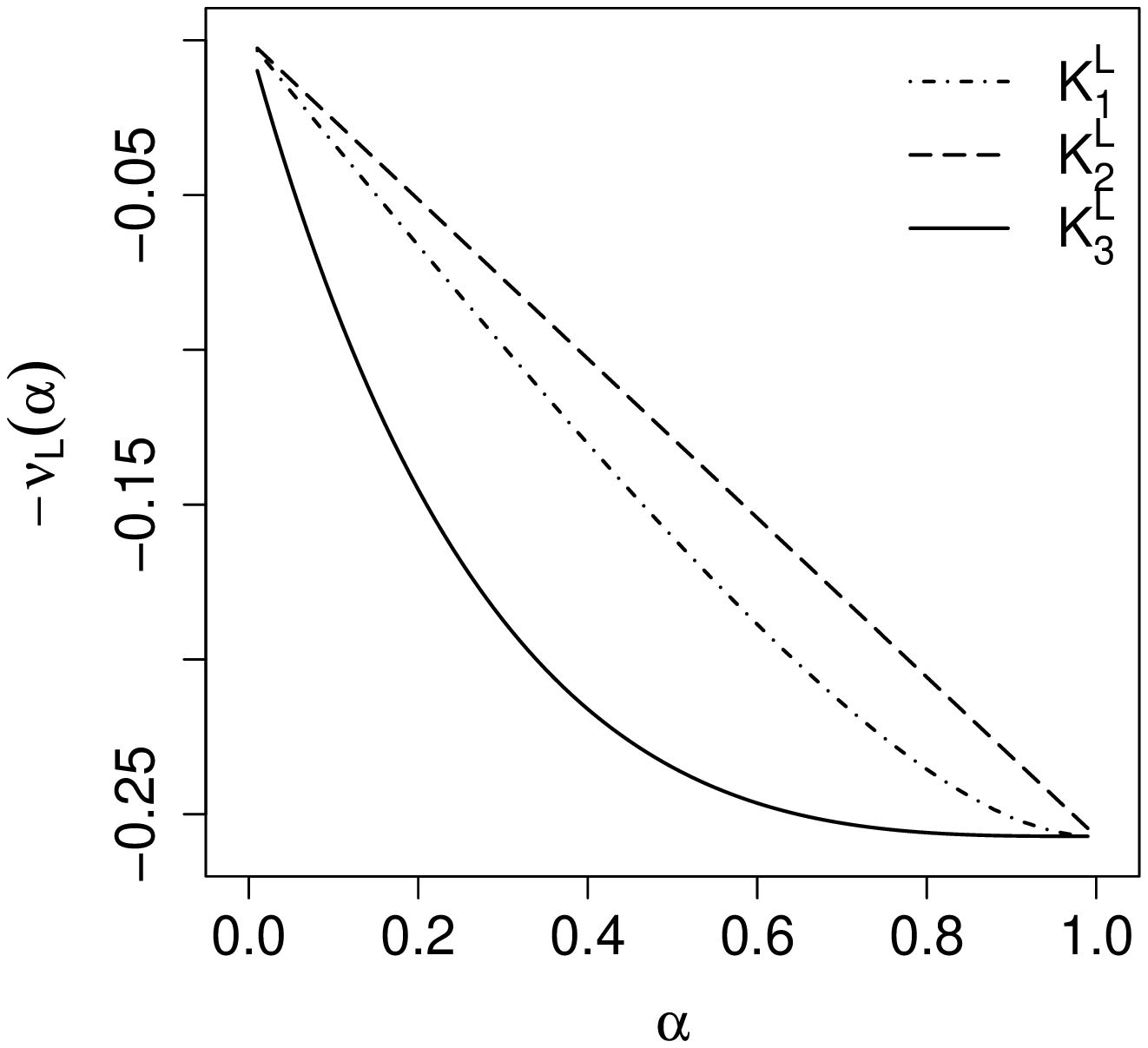} \\*[-1ex]
\end{tabular}

\vspace{-1ex}
\caption{\normalsize\it Functions $\mu_L^2$ and $-\nu_L$ for the left boundary kernels
$K_q^L$, with $q=1,2,3$, where $K$ is the Epanechnikov kernel.} \label{functionsmunu}
\end{figure}

Theorem \ref{theoremboundL} enables us to undertake a first asymptotic comparison between the
boundary kernels $K_q^L$ given by (\ref{K1}), (\ref{K2}) and (\ref{K3}), respectively.
In Figure \ref{functionsmunu} we plot the functions $\mu_L^2$ and $-\nu_L$ which respectively correspond to the coefficients of the most significant
terms in the expansions of the local variance and square bias of estimator $\tF_{nh}(x)$ for $x$ in the left boundary
region. We take for $K$ the Bartlett or Epanechnikov kernel $K(t)=\frac{3}{4}(1-t^2) I(|t|\leq 1)$, but similar conclusions are valid for other polynomial kernels such as the uniform (in this case $K_1^L=K_2^L$),
the biweight or the triweight kernels
\citep[for the definition of these kernels see][p.~31]{WanJ:95}.

From the plots we conclude that the boundary kernel $K^L_3$ has, uniformly over the boundary region,
the biggest asymptotic squared bias but also the lowest asymptotic variance
among the considered boundary kernels. The lowest asymptotic bias is obtained by
$K^L_1$, but this kernel has also the largest asymptotic variance among the considered kernels.
We postpone to Section \ref{compar} the analysis of the combined effect of bias and variance which depends on the underlying
distribution $F$, specially throughout $F^{\prime\prime}(x)^2$ and $F^\prime(x)$ that enter as coefficients of
the terms $\mu_L^2((x-a)/h)$ and $-\nu_L((x-a)/h)$, respectively, in the asymptotic
expansions stated in Theorem \ref{theoremboundL} for the bias and variance of $\tF_{nh}(x)$.

\section{Global behaviour} \label{global}

A widely used measure of the quality of the kernel estimator
is the mean integrated squared error given by
\begin{align*}
\MISE(F;h) & =  \E \int \{ \tF_{nh}(x) - F(x) \}^2  dx
\end{align*}
\begin{align*}
\phantom{\MISE(F;h)} & =  \int \Var \tF_{nh}(x) dx + \int \{ \E \tF_{nh}(x) - F(x) \}^2 dx \\
& =:  \mathbf{V}(F;h) + \mathbf{B}(F;h).
\end{align*}

Next we extend Theorems 4.1 and 4.2 of \citet{Ten:13} by showing that the MISE expansion obtained by
\cite{Jon:90} for the classical kernel estimator (\ref{classical}) is also valid for the
boundary modified kernel estimator (\ref{boundary}) when the left boundary kernel satisfies condition (\ref{C2}).
As before we assume that the right boundary kernel $K^R$ is given by $K^R(u;\alpha)=K^L(-u;\alpha)$, for $u\in\R$ and
$\alpha\in\,]0,1[$.

\medskip

\begin{theorem} \label{misetheo}
If $K^L(u;\alpha)$ satisfies condition (\ref{C2}) with
\begin{align} \label{ker3F}
\int_0^1  |\mu_{0,L}|(\alpha)^2 d\alpha < \infty,
\end{align}
and the restriction of $F$ to the interval $[a,b]$ is twice continuously differentiable, we have:
$$\mathbf{V}(F;h)= \frac{1}{n} \int F(x)(1-F(x)) dx - \frac{h}{n} \int uB(u) du + O\left( n^{-1} h^2 \right)$$
and
$$\mathbf{B}(F;h)= \frac{h^4}{4}\, \left( \int u^2 K(u) du \right)^2 
\int F^{\prime\prime}(x)^2 dx + o\left(h^4 \right).$$
Moreover, if $F$ is not the uniform distribution function on $[a,b]$,
the asymptotically optimal bandwidth, in the sense of minimising the MISE expansion leading terms, is given by
\begin{align*}
h_0= \delta(K) \left( \int F^{\prime\prime}(x)^2 dx \right)^{-1/3} n^{-1/3},
\end{align*}
where
$$\delta(K)=\left(  \int uB(u)\, du \right)^{1/3} \left( \int u^2 K(u) du \right)^{-2/3}.$$
\end{theorem}

\medskip

A classical measure of a distribution function estimator performance is the supremum distance between
such an estimator and the underlying distribution function $F$.
Next we extend Theorems 3.1 and 3.2 of \citet{Ten:13} by establishing
the almost complete uniform convergence and the Chung-Smirnov law
of iterated logarithm for kernel estimator (\ref{boundary}).
These properties have been first obtained for estimator (\ref{classical}) by
\citet{Nad:64}, \citet{Win:73,Win:79} and \citet{Yam:73}.
We denote by $||\cdot||$ the supremum norm.

\medskip

\begin{theorem} \label{supremum}
If $K^L(u;\alpha)$ is such that
$$
\sup_{\alpha\in\, ]0,1[} |\mu_{0,L}|(\alpha) < \infty,
$$
we have
$$||\tF_{nh} - F || \rightarrow 0 \quad \mbox{almost completely}.$$
Additionally, if $F$ is Lipschitz on $[a,b]$ and $(n/\log\log n)^{1/2} h \rightarrow 0$,
then $\tF_{nh}$ has the Chung-Smirnov property, i.e.,
$$\limsup_{n\rightarrow\infty} \, (2n/\log\log n)^{1/2} || \tF_{nh} - F || \leq 1 \;\mbox{ almost surely}.$$
The same is true under the less restrictive condition $(n/\log\log n)^{1/2} h^2 \rightarrow 0$,
whenever $K^L$ satisfies (\ref{C2}) and $F^\prime$ is Lipschitz on $[a,b]$.
\end{theorem}

\begin{remark}
{\rm
The asymptotically optimal bandwidth $h_0$ given in Theorem \ref{misetheo} satisfies condition
$(n/\log\log n)^{1/2} h^2 \rightarrow 0$, but not condition $(n/\log\log n)^{1/2} h \rightarrow 0$.
}
\end{remark}

\section{Exact finite sample comparisons} \label{compar}

\begin{figure}[!p]

\centering
\begin{tabular}{c@{\hspace{0.5em}}c@{\hspace{1.5em}}c}
 & \hspace{5mm}$F^{\prime\prime}_+(0)=6$ & \hspace{5mm}$F^{\prime\prime}_+(0)=30$  \\*[0.5ex]
\rotatebox{90}{\hspace{17mm}$F_+^\prime(0)=0$} & \includegraphics[scale=0.375]{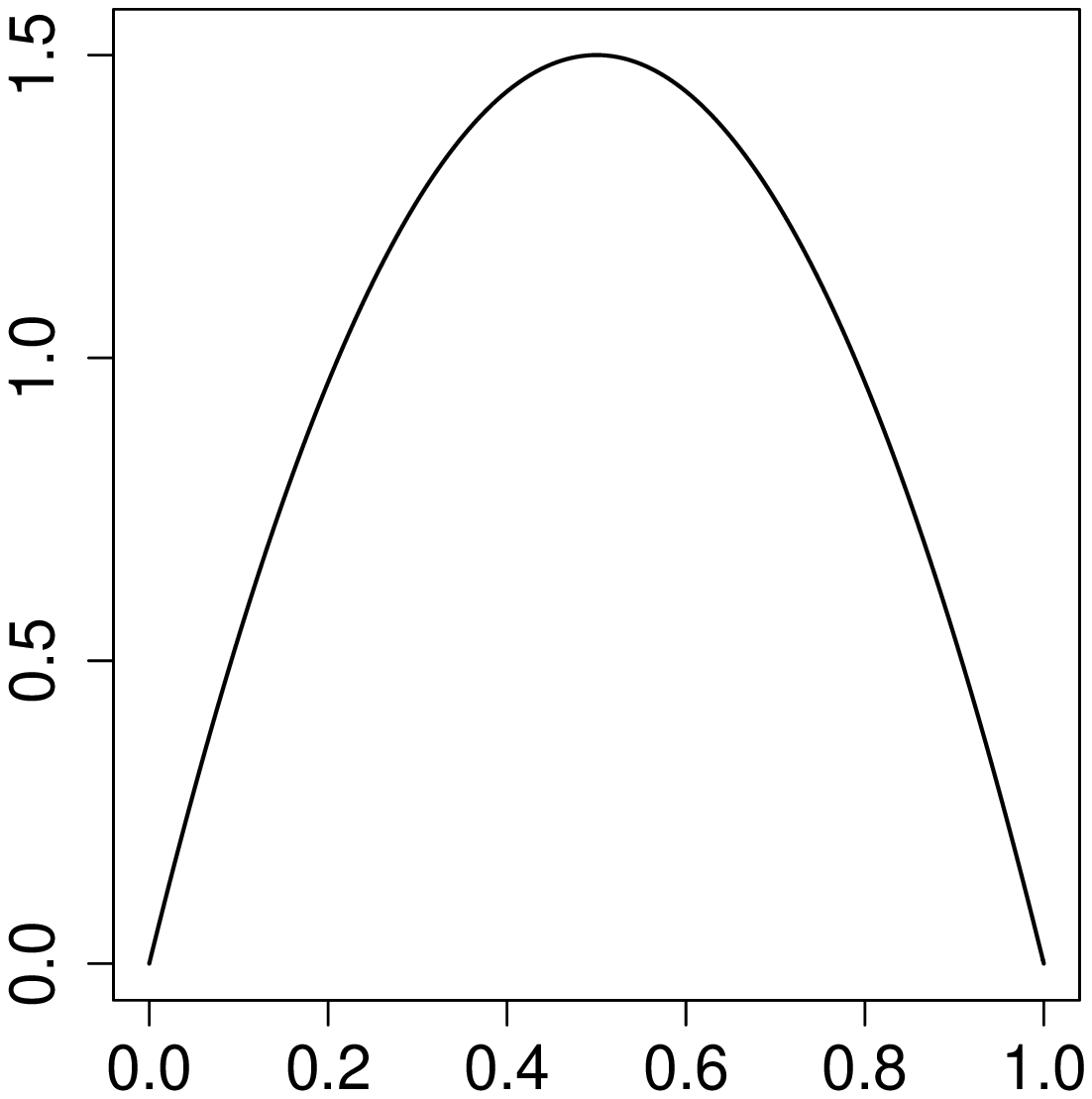} & \includegraphics[scale=0.375]{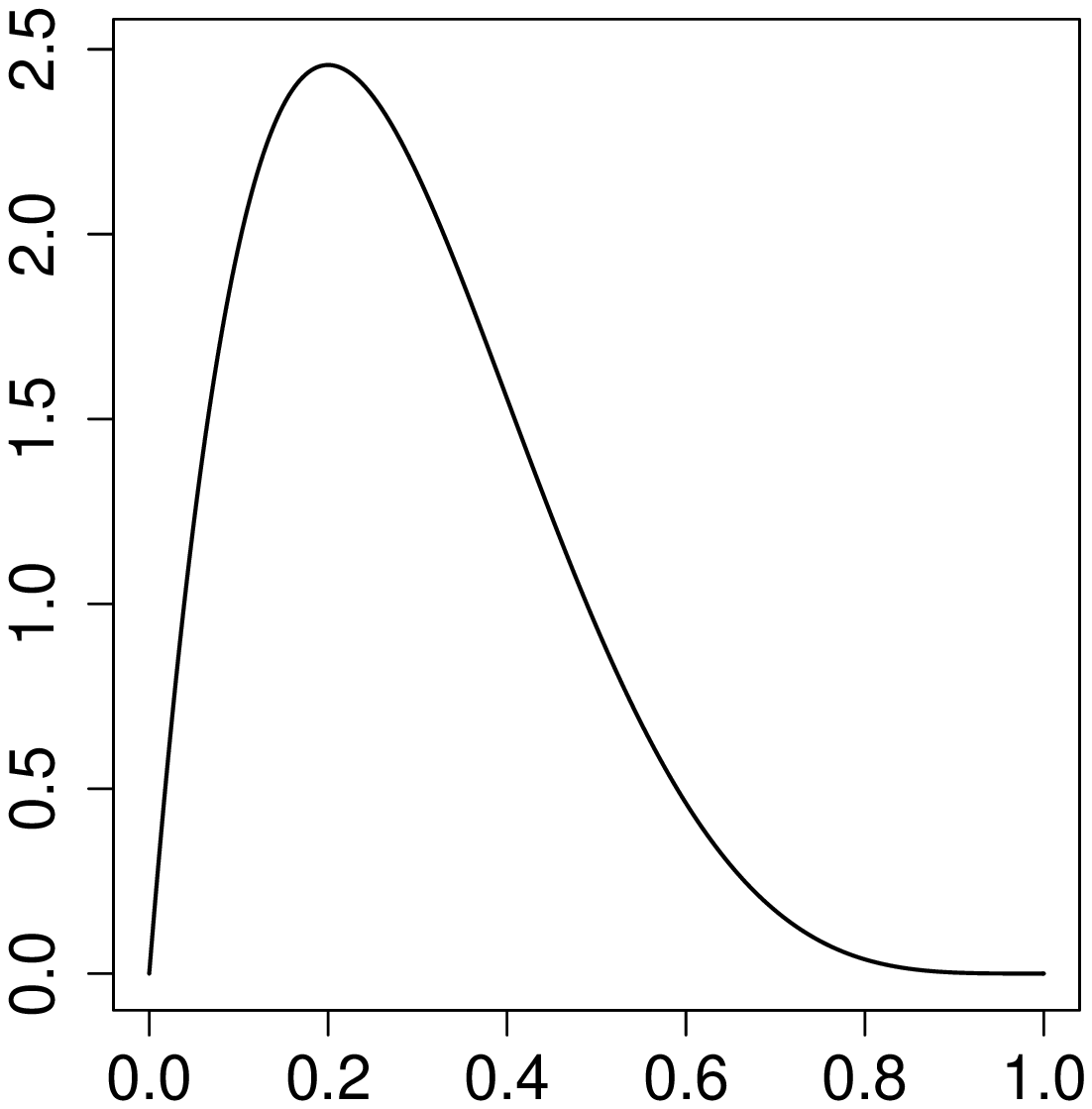} \\*[-1.5ex]
\rotatebox{90}{\hspace{17mm}$F_+^\prime(0)=0.5$} & \includegraphics[scale=0.375]{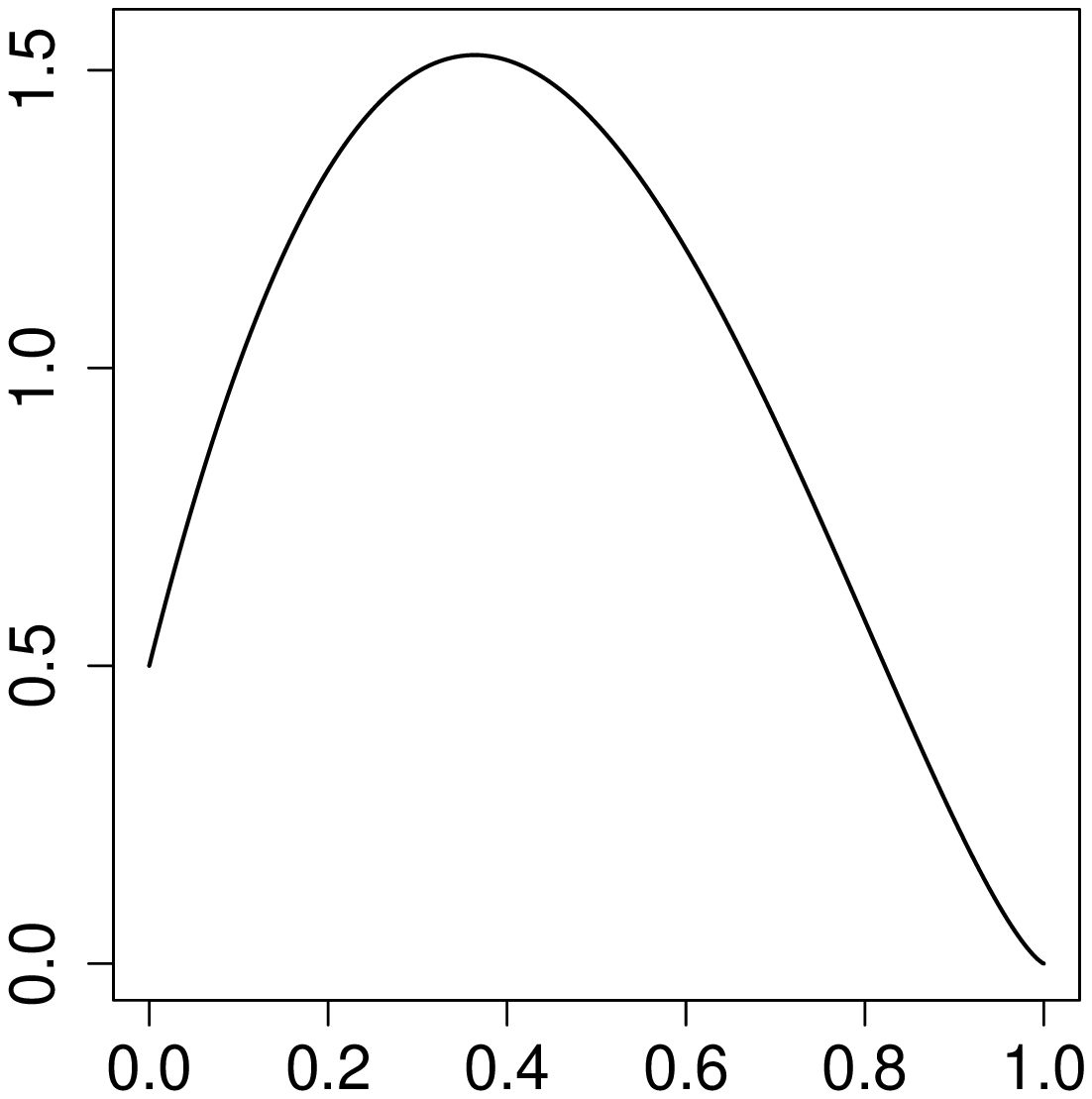} & \includegraphics[scale=0.375]{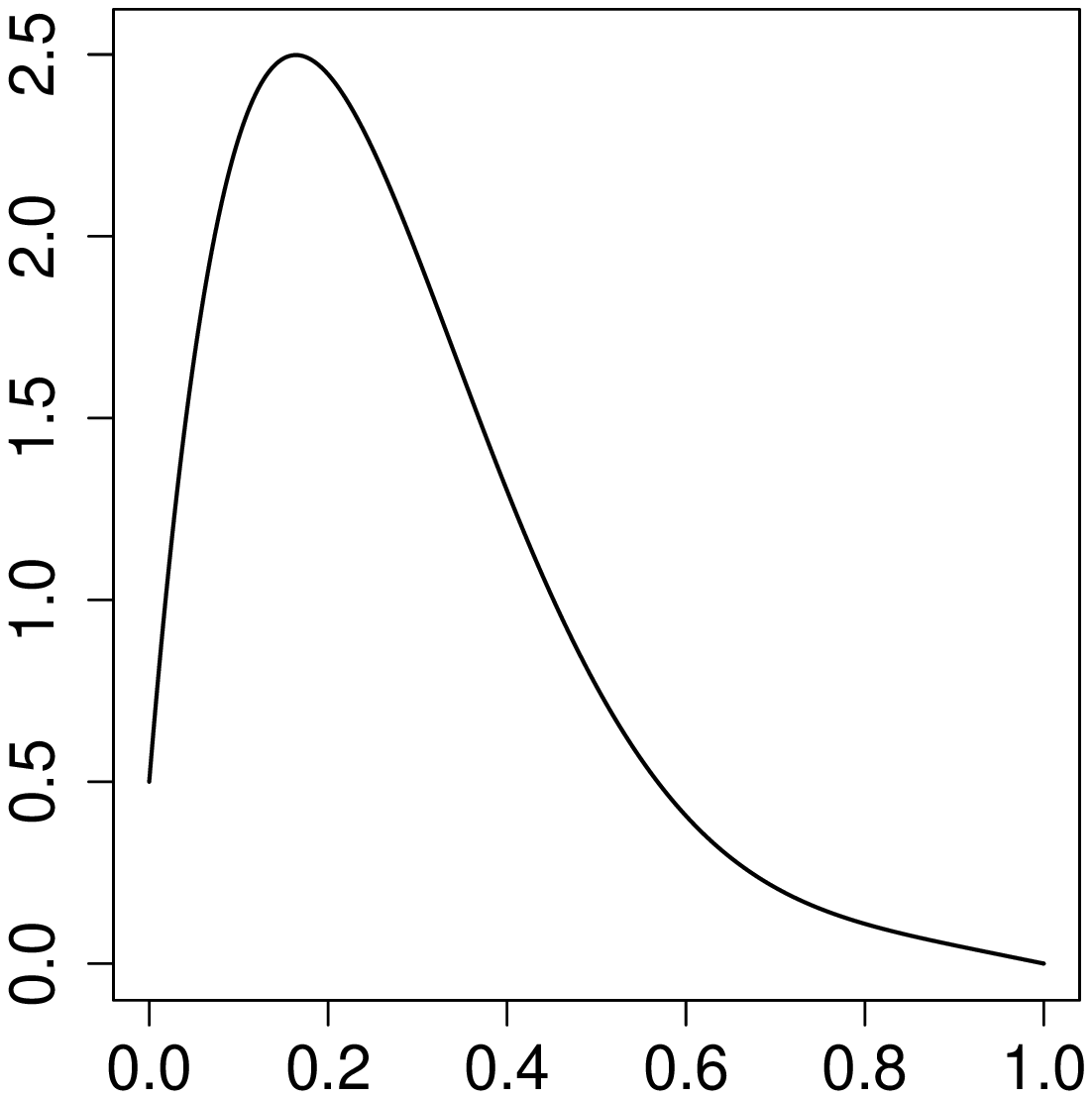} \\*[-1.5ex]
\rotatebox{90}{\hspace{17mm}$F_+^\prime(0)=1$} & \includegraphics[scale=0.375]{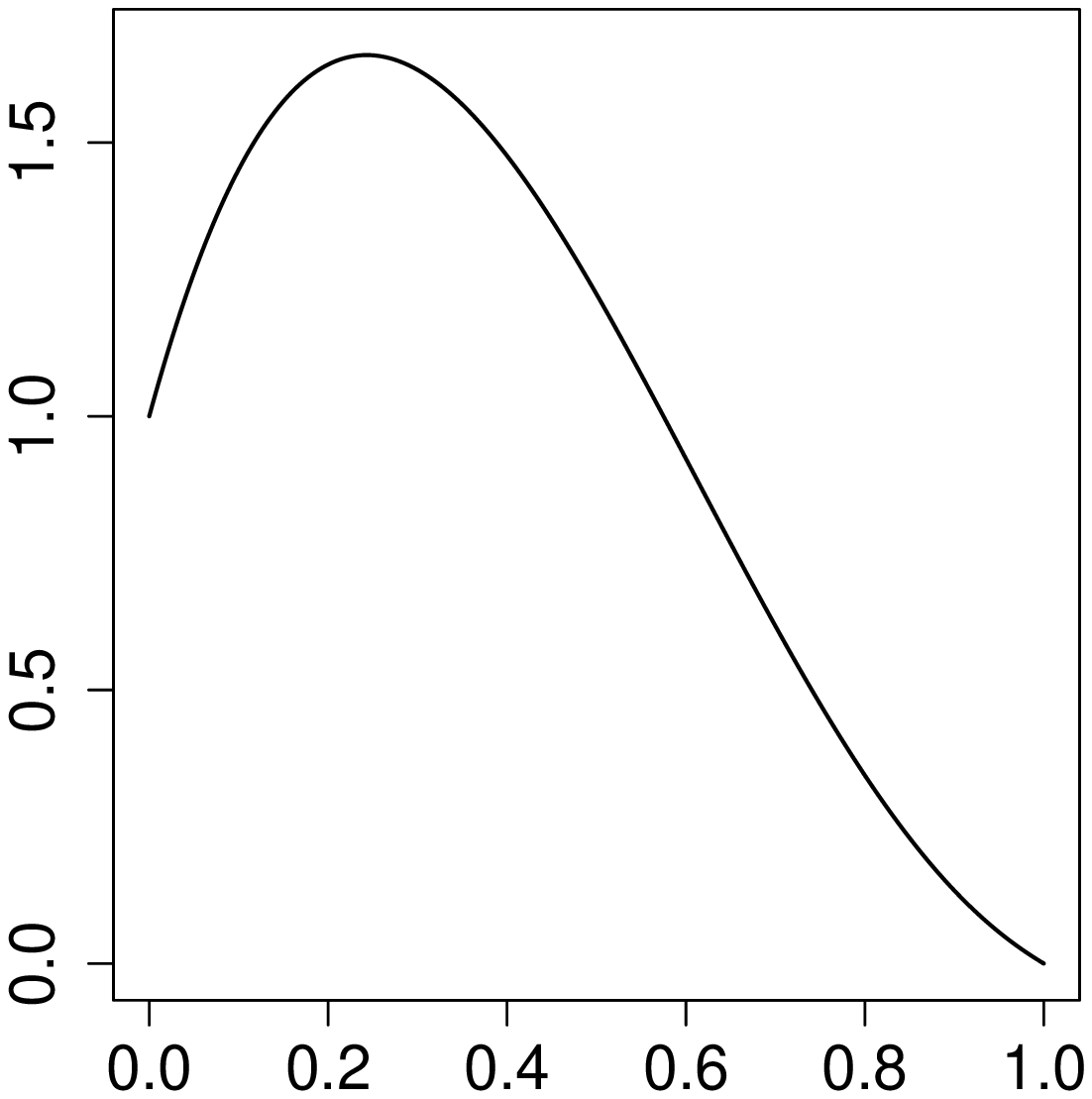} & \includegraphics[scale=0.375]{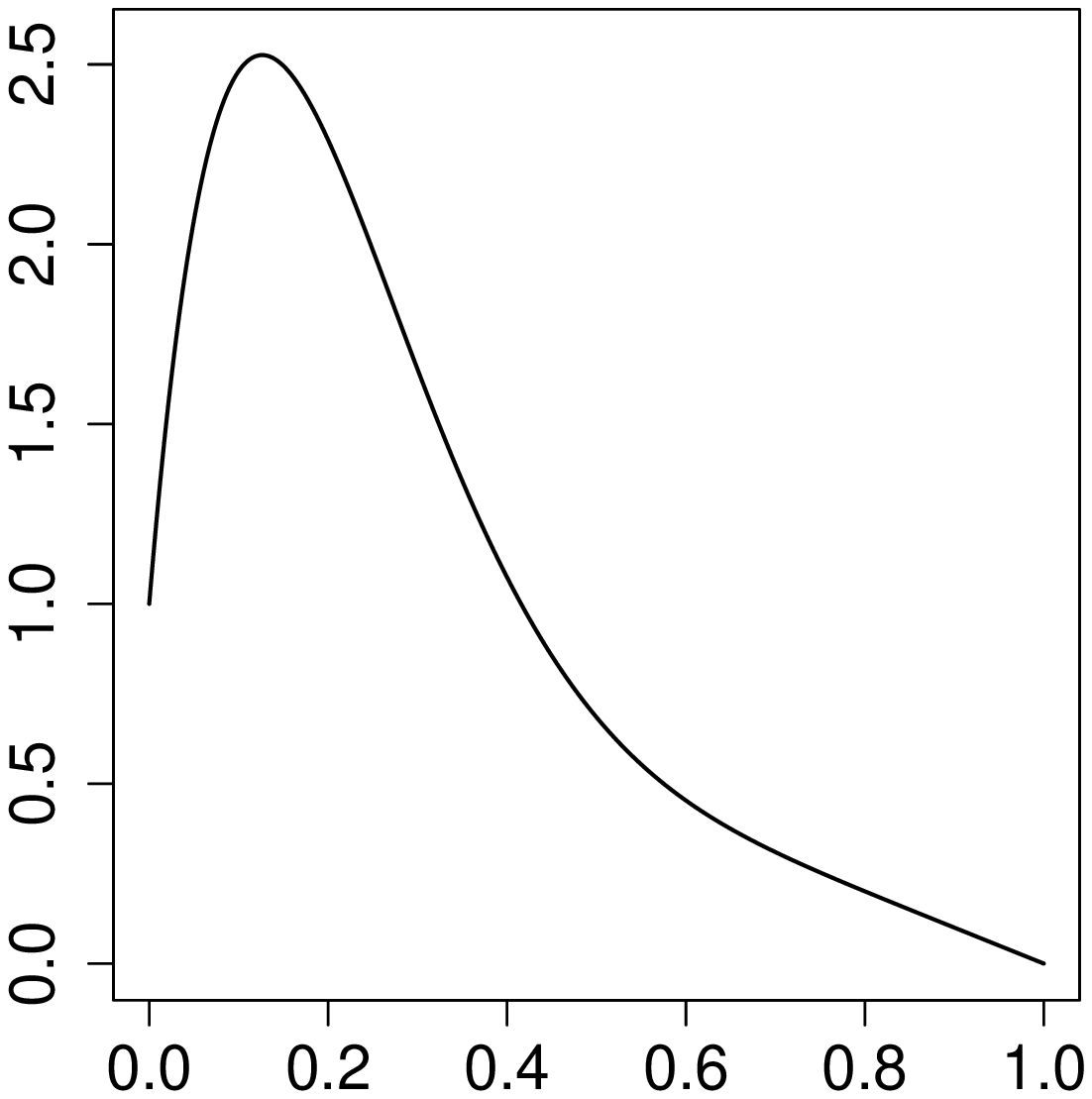} \\*[-1.5ex]
\rotatebox{90}{\hspace{17mm}$F_+^\prime(0)=1.5$} & \includegraphics[scale=0.375]{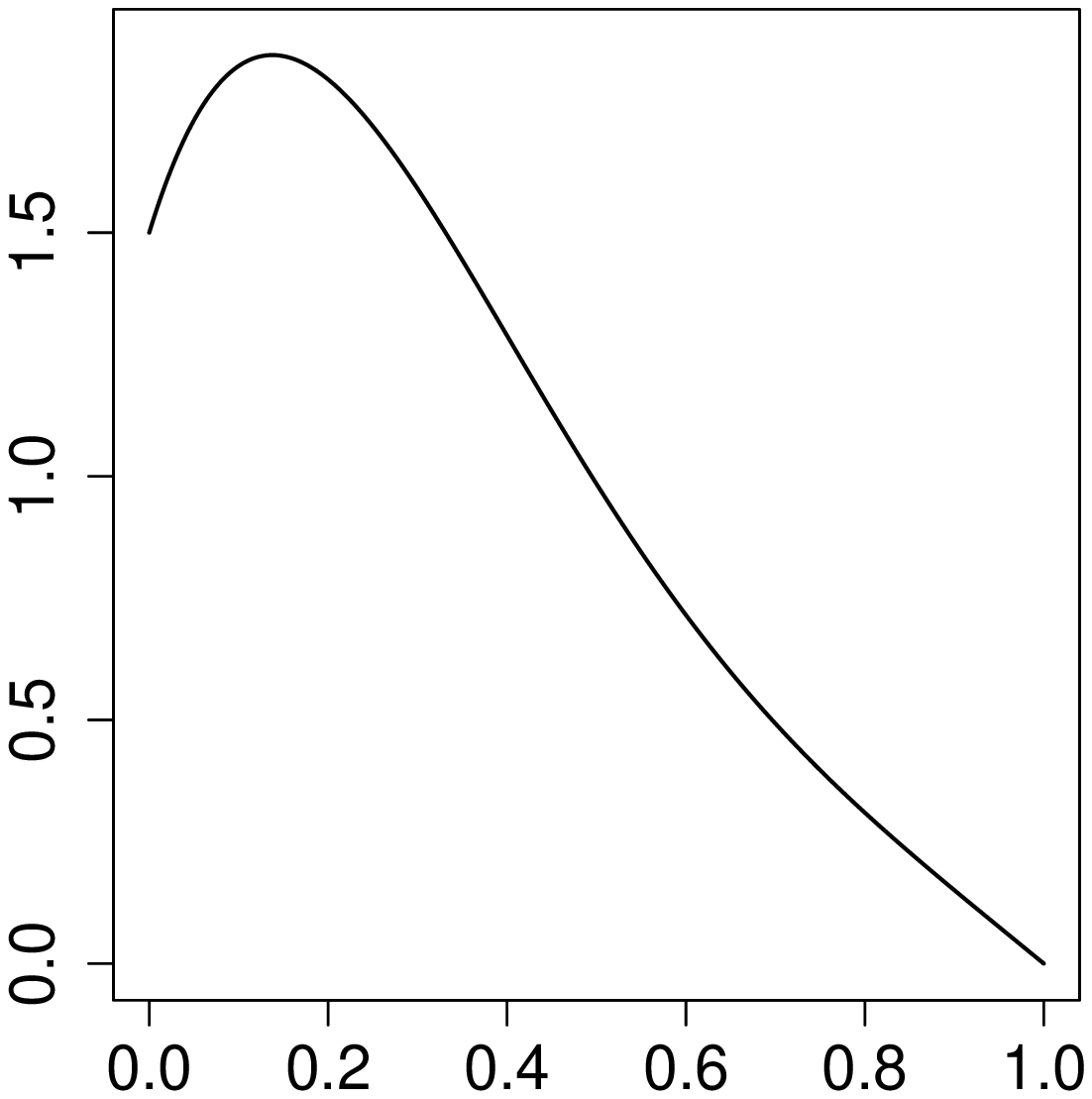} & \includegraphics[scale=0.375]{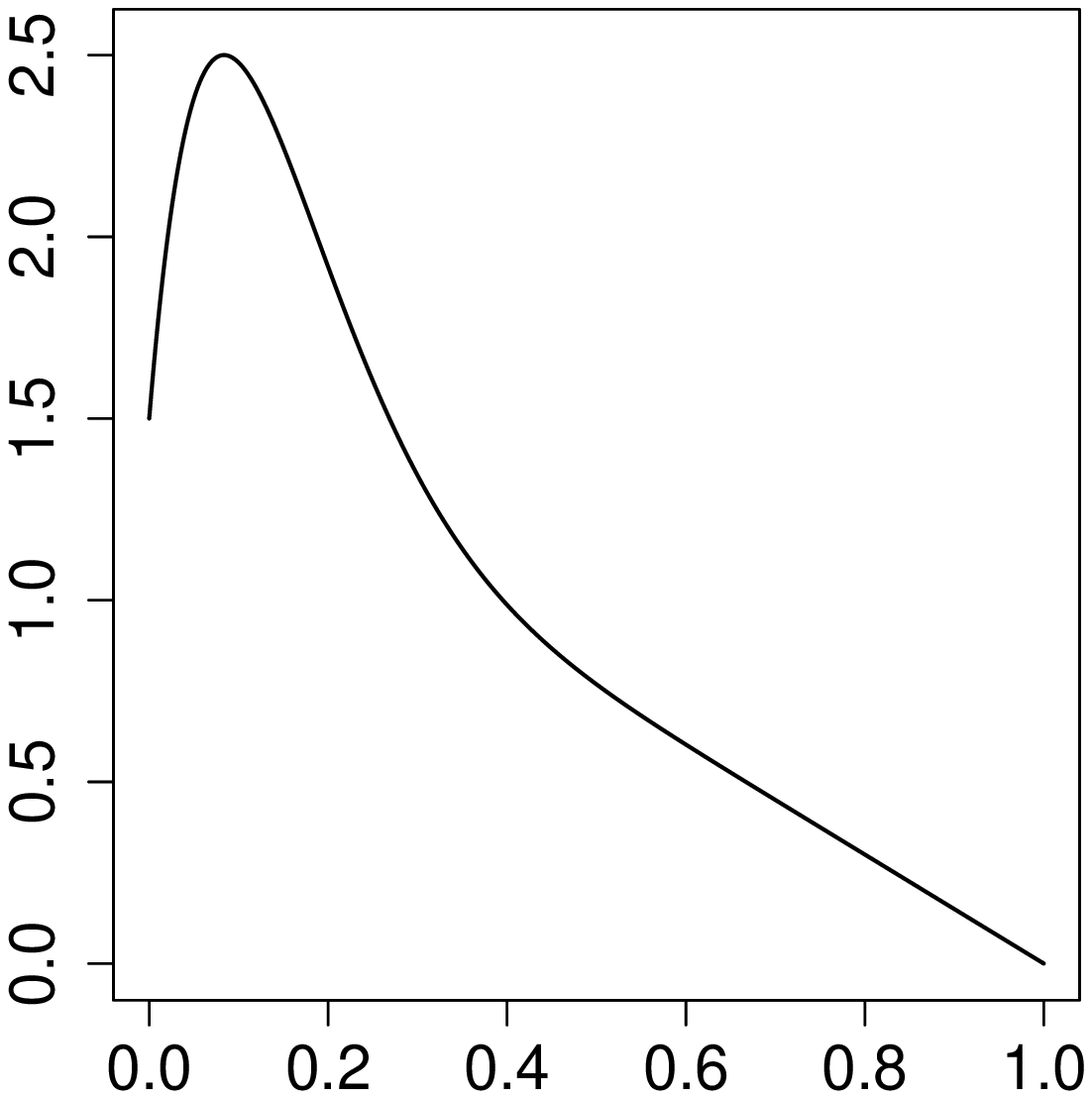}
\end{tabular}

\vspace{-1ex}
\caption{\normalsize\it Beta mixture densities $w B(1,2)+(1-w) B(2,b)$ with
$F_+^\prime(0)=0, 0.5, 1, 1.5$ and
$F^{\prime\prime}_+(0)=6$ (left column) and $F^{\prime\prime}_+(0)=30$ (right column).} \label{testdist}
\end{figure}

In this section we compare the boundary performance of the
kernel estimator $\tF_{nh}$ when we take for $K^L$ one of the left boundary kernels
given by (\ref{K1}), (\ref{K2}) and (\ref{K3}), respectively.
For that, we have used as test distributions some beta mixtures of the form
$w B(1,2)+(1-w) B(2,b),$ where $w\in\,[0,1]$ and the shape parameter $b$ is such that $b\geq 2$.
Four values of $w=0,0.25,0.5,0.75$ were considered, which
lead to distributions with $F_+^\prime(0)=0,0.5,1,1.5$, respectively.
For each one of the previous weights $w$, two values for the shape parameter $b$ were taken
in order to get a second order derivative $F^{\prime\prime}_+(0)$ equal to $6$ and $30$.
The considered set of test distributions is shown in Figure \ref{testdist}.

\begin{figure}[!p]

\centering
\begin{tabular}{c@{\hspace{1.5em}}c@{\hspace{2em}}c}
 & \hspace{5mm}$F^{\prime\prime}_+(0)=6$ & \hspace{5mm}$F^{\prime\prime}_+(0)=30$  \\*[0.5ex]
\rotatebox{90}{\hspace{17mm}$F_+^\prime(0)=0$} & \includegraphics[scale=0.375]{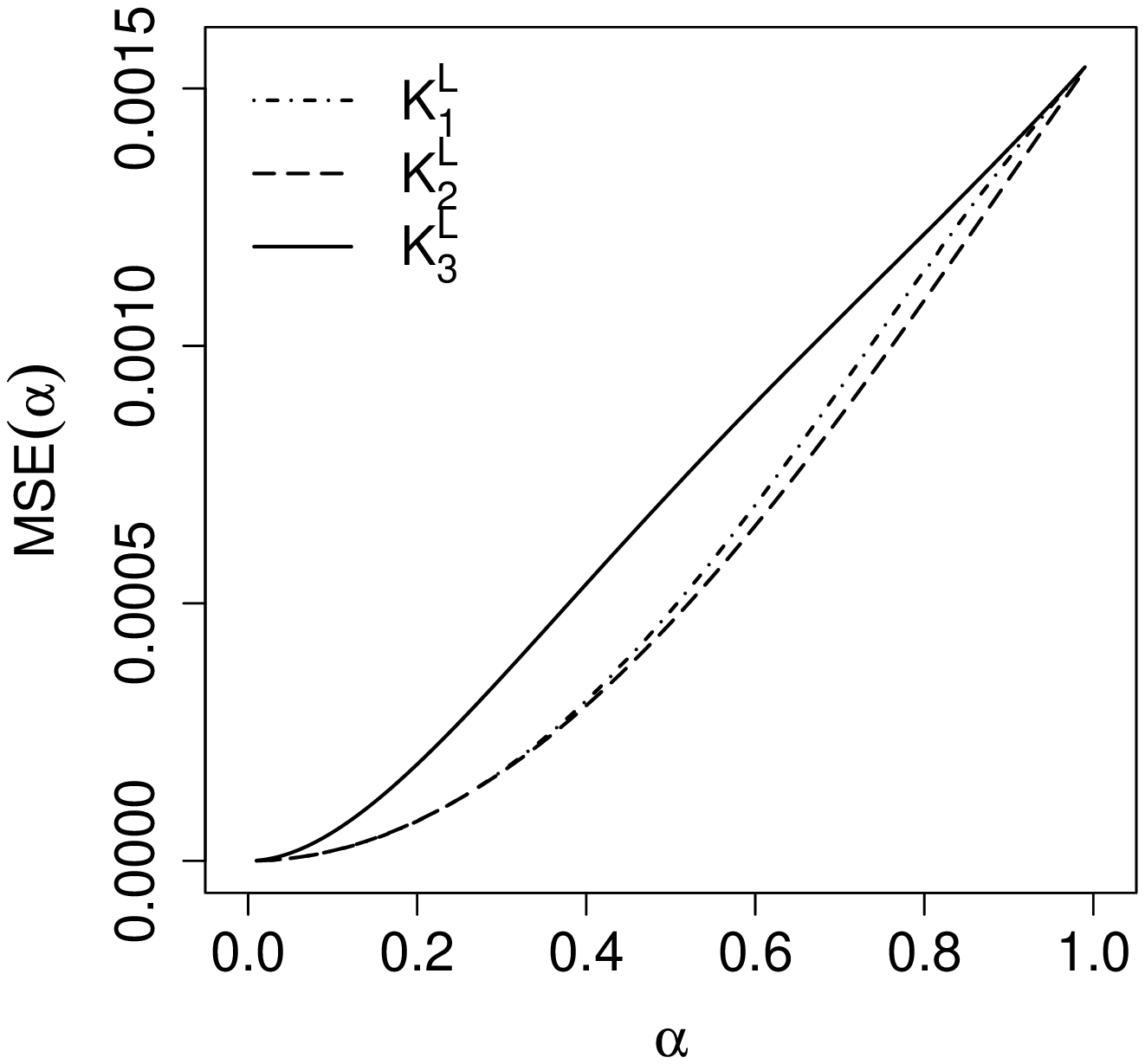} & \includegraphics[scale=0.375]{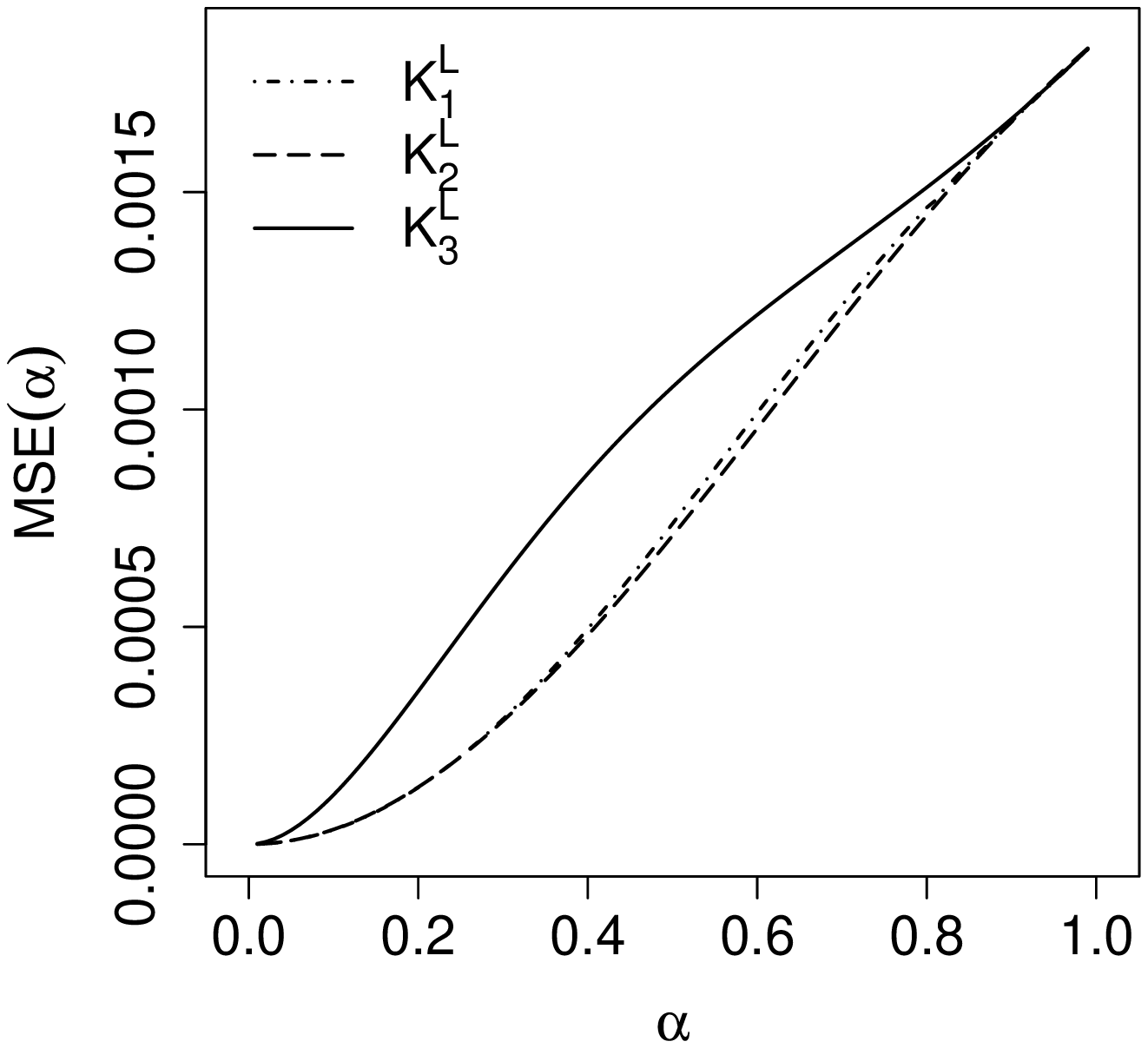} \\*[-0.5ex]
\rotatebox{90}{\hspace{17mm}$F_+^\prime(0)=0.5$} & \includegraphics[scale=0.375]{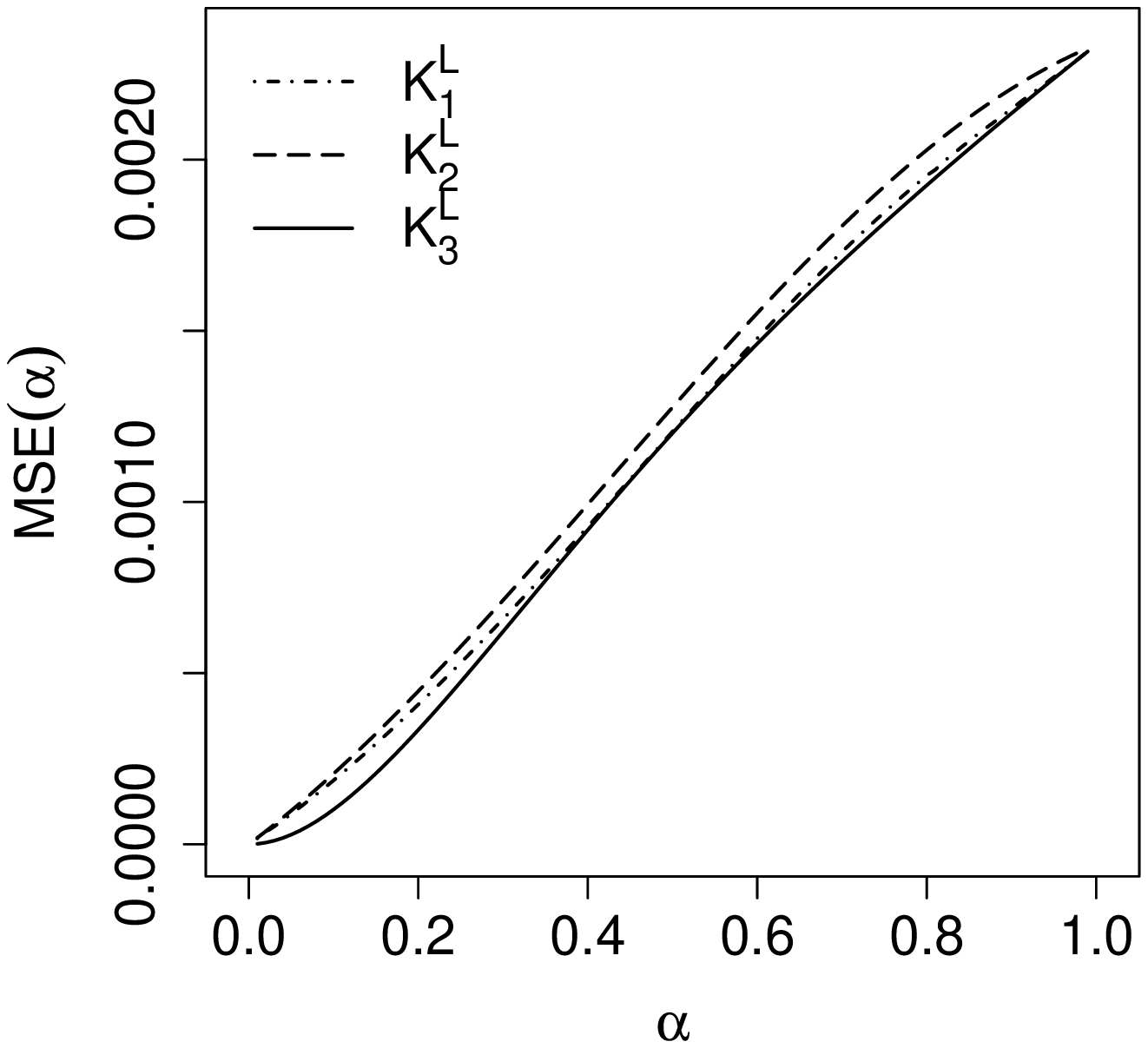} & \includegraphics[scale=0.375]{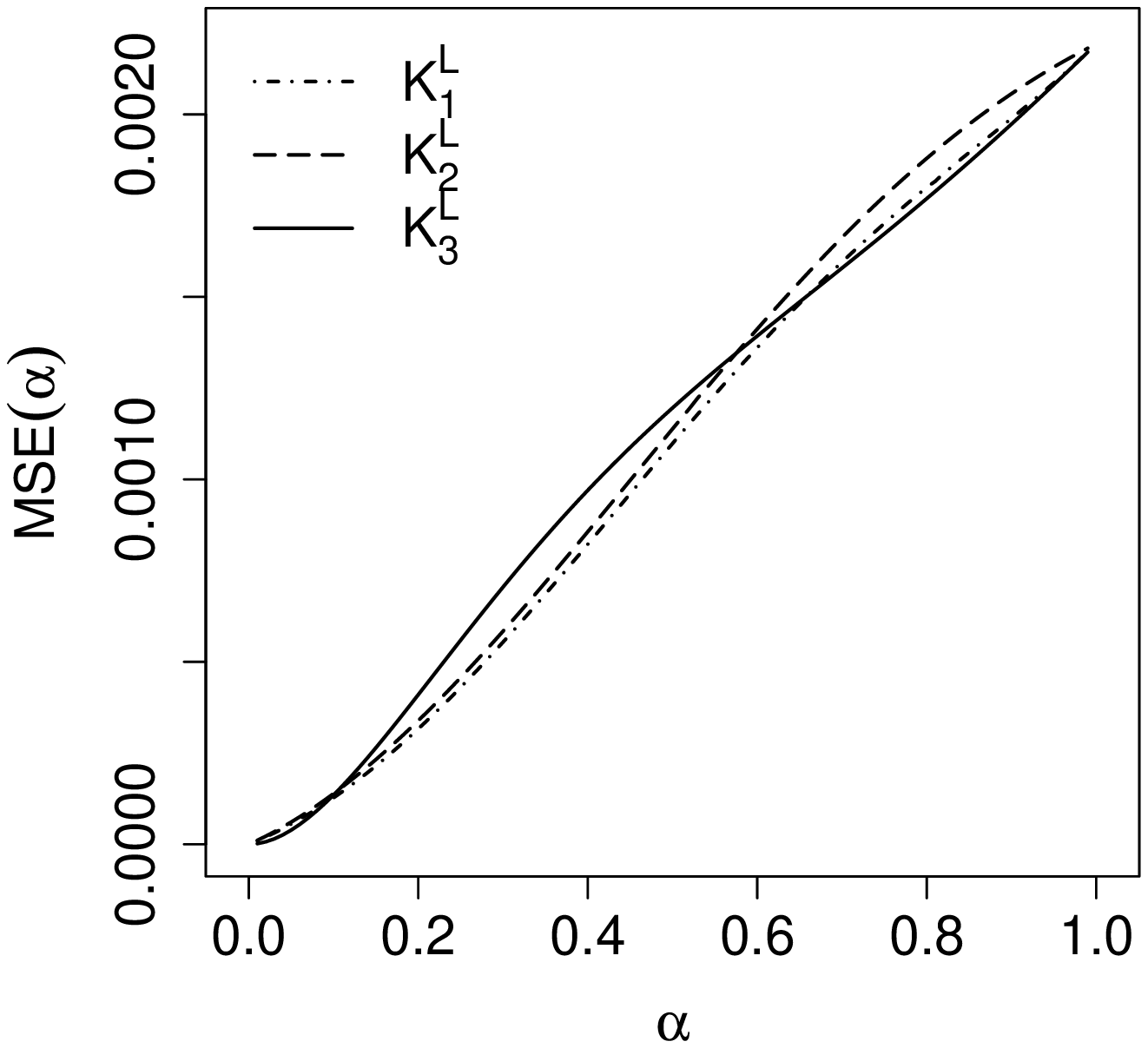} \\*[-0.5ex]
\rotatebox{90}{\hspace{17mm}$F_+^\prime(0)=1$} & \includegraphics[scale=0.375]{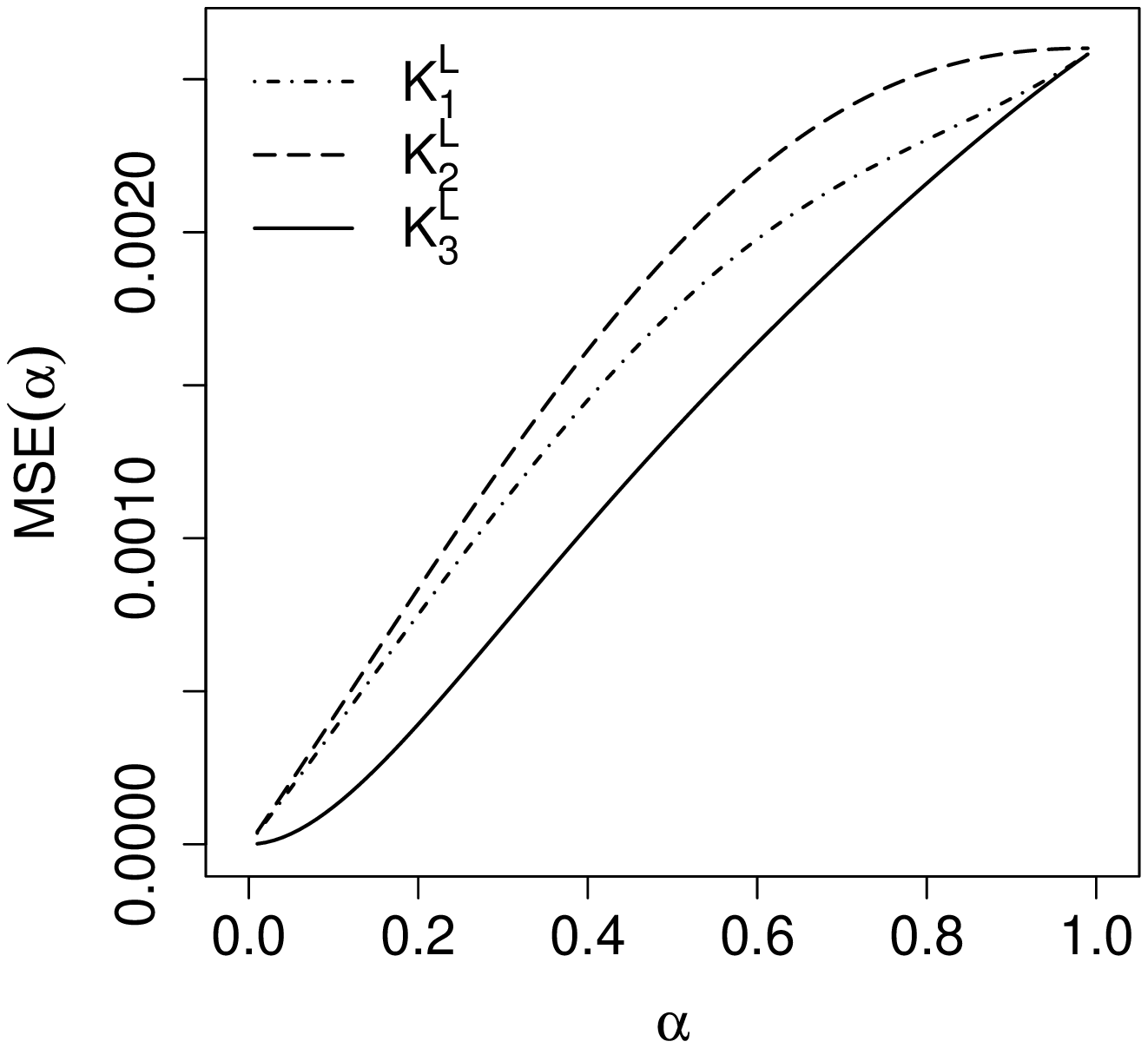} & \includegraphics[scale=0.375]{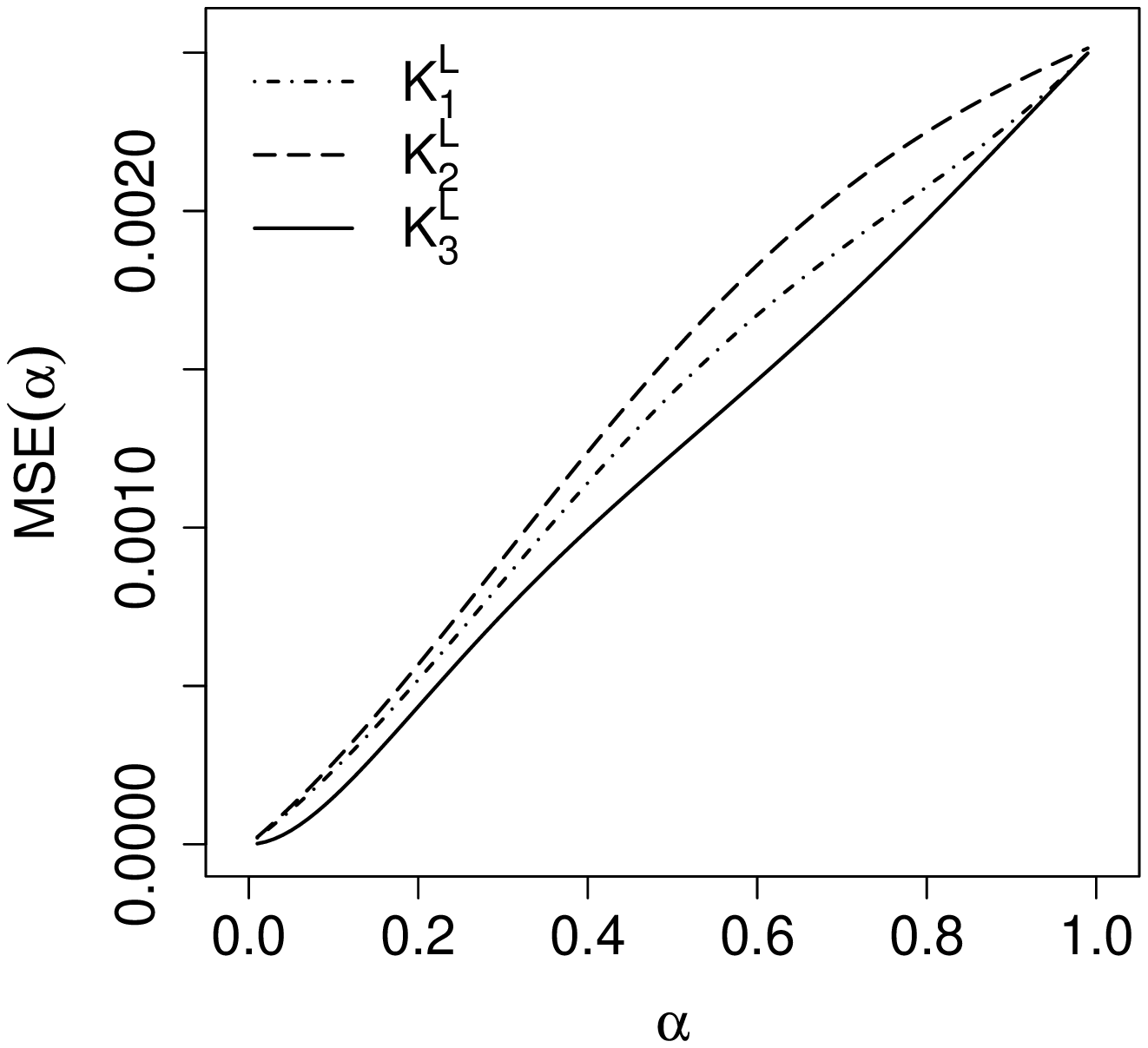} \\*[-0.5ex]
\rotatebox{90}{\hspace{17mm}$F_+^\prime(0)=1.5$} & \includegraphics[scale=0.375]{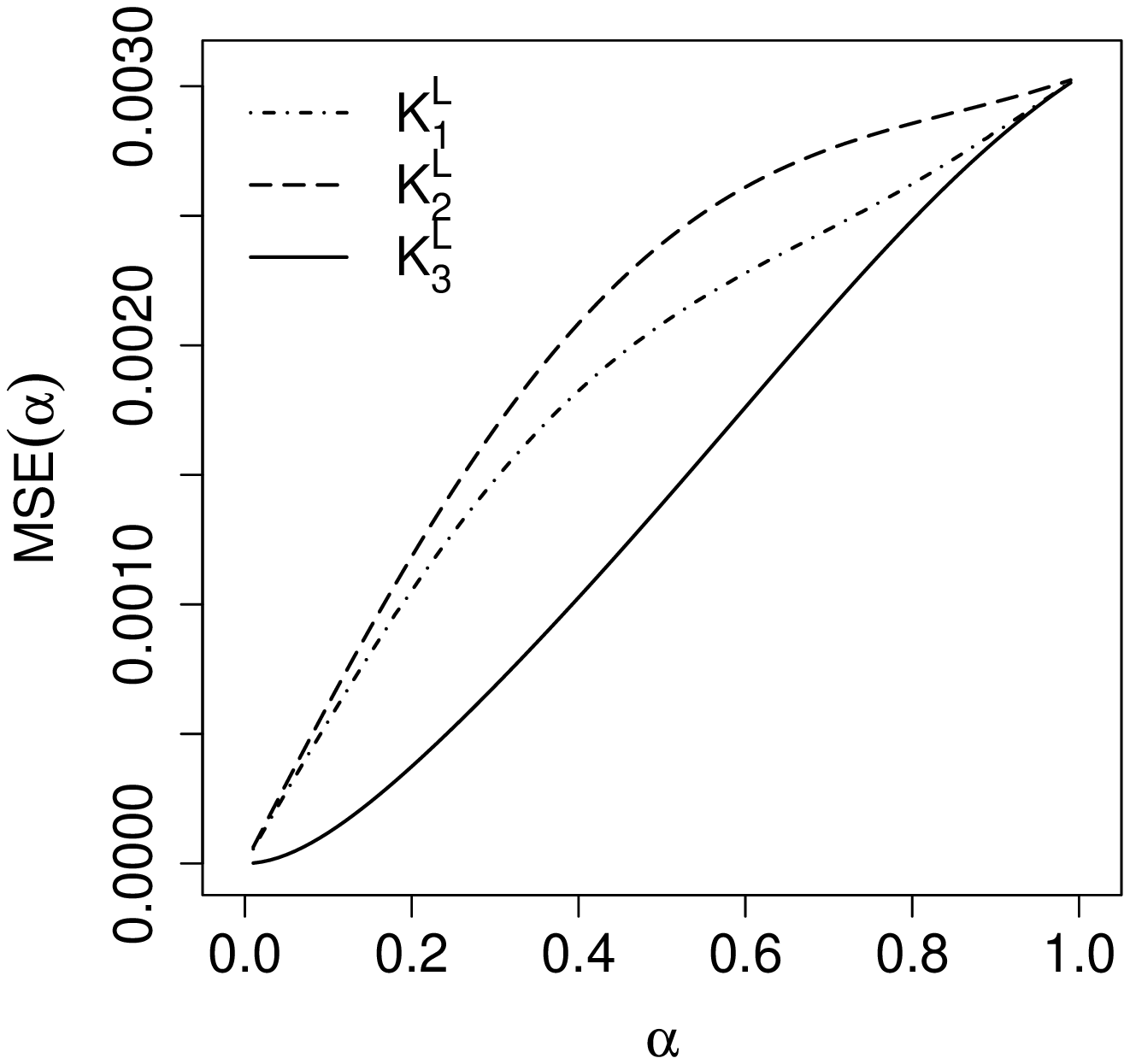} & \includegraphics[scale=0.375]{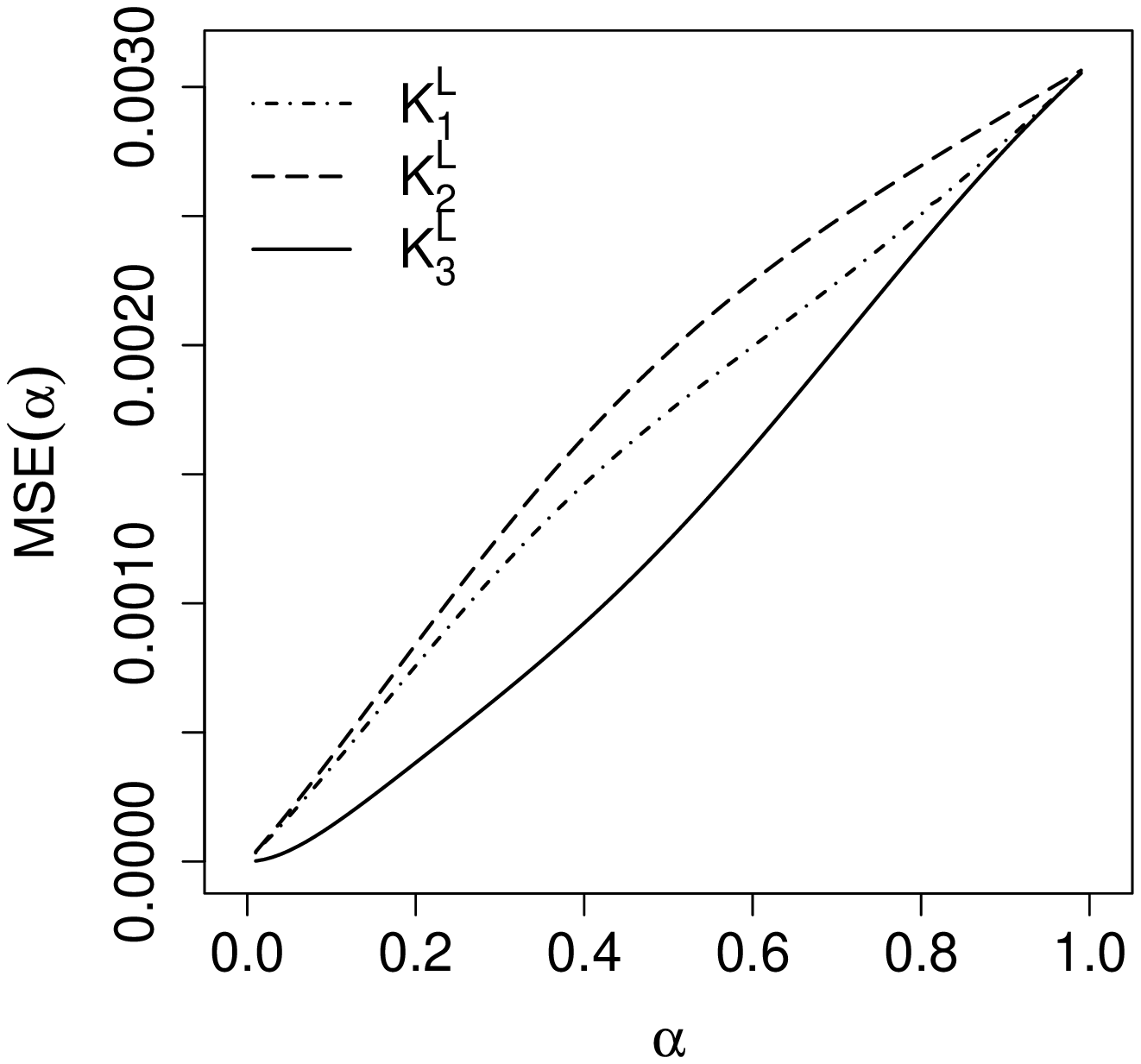}
\end{tabular}

\vspace{-1ex}
\caption{\normalsize\it $\mathrm{MSE}(\alpha)$ for $K^L_q$, $q=1,2,3$, with $K$ the Epanechnikov kernel,
where $F$ is the beta mixture distribution $w B(1,2)+(1-w) B(2,b)$ with
$F^{\prime}_+(0)=0, 0.5, 1, 1.5$, $F^{\prime\prime}_+(0)=6$ (left column) and $F^{\prime\prime}_+(0)=30$ (right column). The sample size is $n=50$.} \label{mseboundary}
\end{figure}

For each one of these test distributions we present in Figure \ref{mseboundary}
the exact mean square error of $\tF_{nh}(x)$, for
$x=\alpha h$ and $\alpha\in\,]0,1[$, given by
$$\mathrm{MSE}(\alpha)=\mathrm{V}(\alpha) + \mathrm{B}(\alpha)^2,$$
where
\begin{align*}
n \mathrm{V}(\alpha) & := n \Var \tF_{nh}(a+\alpha h)
= \! \int F(a+(\alpha-u)h)B^L(u;\alpha) du - \big( \E \tF_{nh}(a+\alpha h) \big)^2
\end{align*}
and
\begin{align*}
\mathrm{B}(\alpha)
& := \E \tF_{nh}(a+\alpha h) - F(a+\alpha h)
= \! \int F(a+(\alpha-u)h) K^L(u;\alpha) \,du - F(a+\alpha h)
\end{align*}
(on these expressions see Section \ref{proofs} below).
The global bandwidth $h$ that determines the boundary region was always taken equal to
the asymptotically optimal bandwidth $h_0$ given in Theorem \ref{misetheo}, and
we have considered the sample size $n=50$.
Similar pictures were generated for sample sizes $n=100$ and $n=200$,
but they were not included here to save space.
As before, we have taken for $K$ the Epanechnikov kernel.

From the graphics we conclude that the boundary behaviour of
the kernel estimator based on the boundary kernels
$K^L_q$, for $q=1,2,3$, is dominated by the
magnitude of the underlying density $f=F^\prime$ over the boundary region.
For large values of $F_+^\prime(0)$ we see that
the boundary kernel $K_3^L$ is superior to both $K_1^L$ and $K_2^L$, being the advantage over
the second order boundary kernels bigger for large than for small values of $F_+^{\prime\prime}(0)^2$.
Notice that this latter conclusion is in accordance with the asymptotic
comparisons presented in Section \ref{secboundary}.
Although less performing than $K_3^L$, the kernel $K_1^L$ is, in this case, superior to $K_2^L$.
When the underlying density is such that $F_+^\prime(0)=0$,
in which case the classical kernel estimator does not suffer from boundary problems, we see that
the boundary kernels $K_1^L$ and $K_2^L$ perform similarly being both slightly better than $K_3^L$.
Finally, for intermediate values of $F_+^\prime(0)$
the three considered left boundary kernels are equally performing.
Based on this analysis, we conclude that none of the considered
boundary kernels is the best over the considered
set of test distributions. However, the kernel $K_3^L$ shows to be particularly
interesting because it is especially performing when the classical boundary
kernel estimator suffers from severe boundary problems.

\begin{figure}[!t]

\centering
\begin{tabular}{c@{\hspace{.1em}}c@{\hspace{.1em}}c}
\includegraphics[scale=0.375]{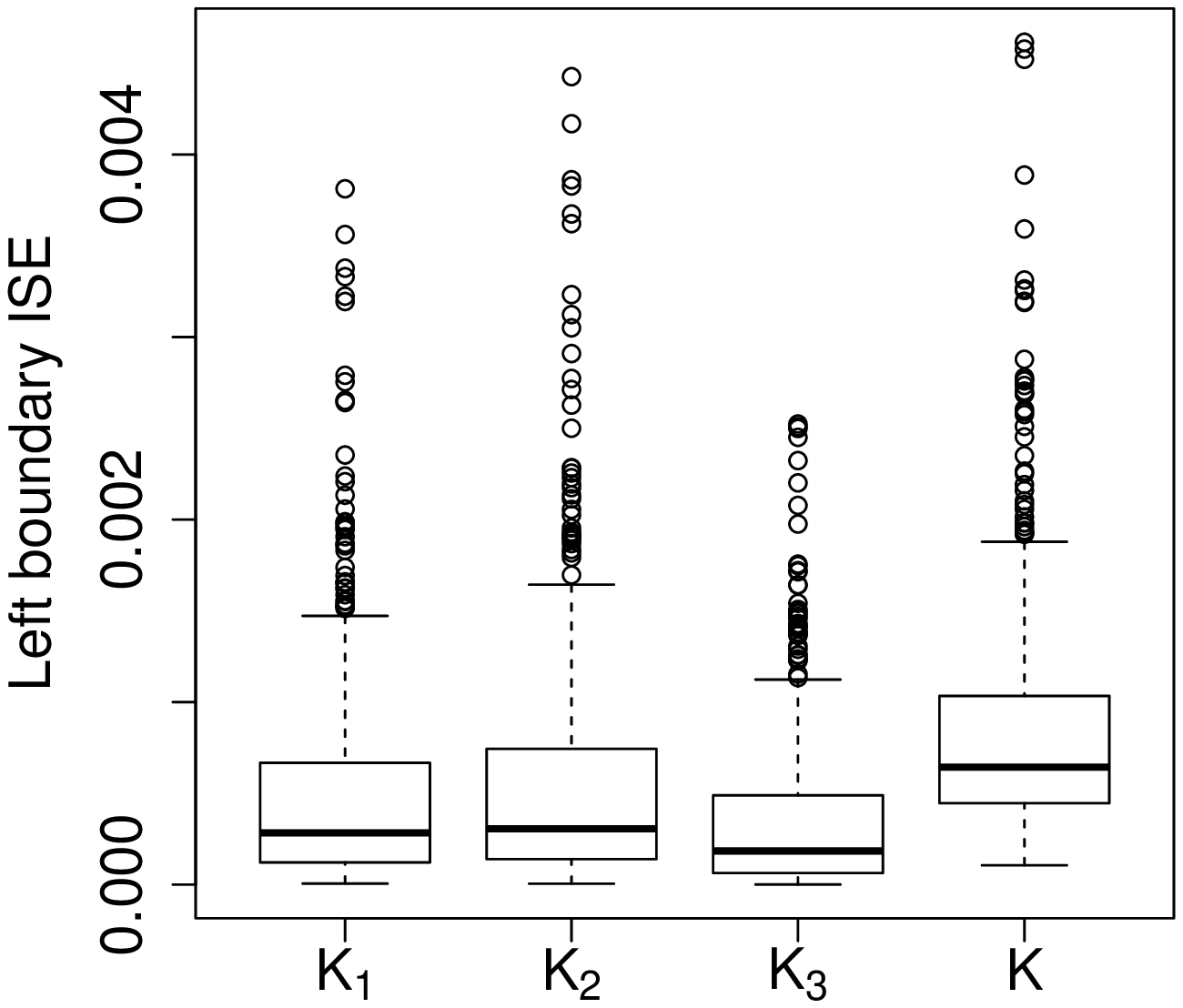} & \includegraphics[scale=0.375]{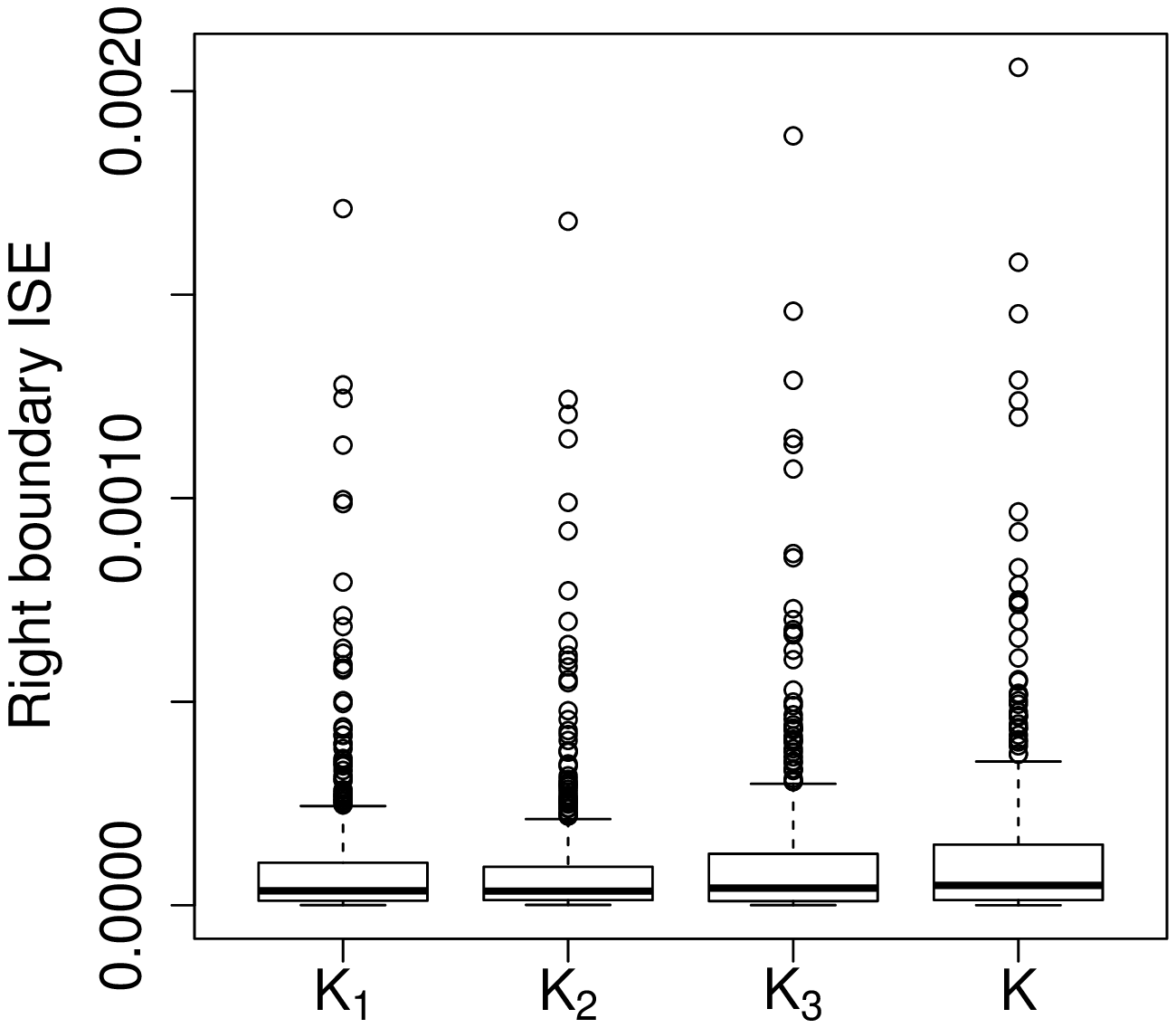} & \includegraphics[scale=0.375]{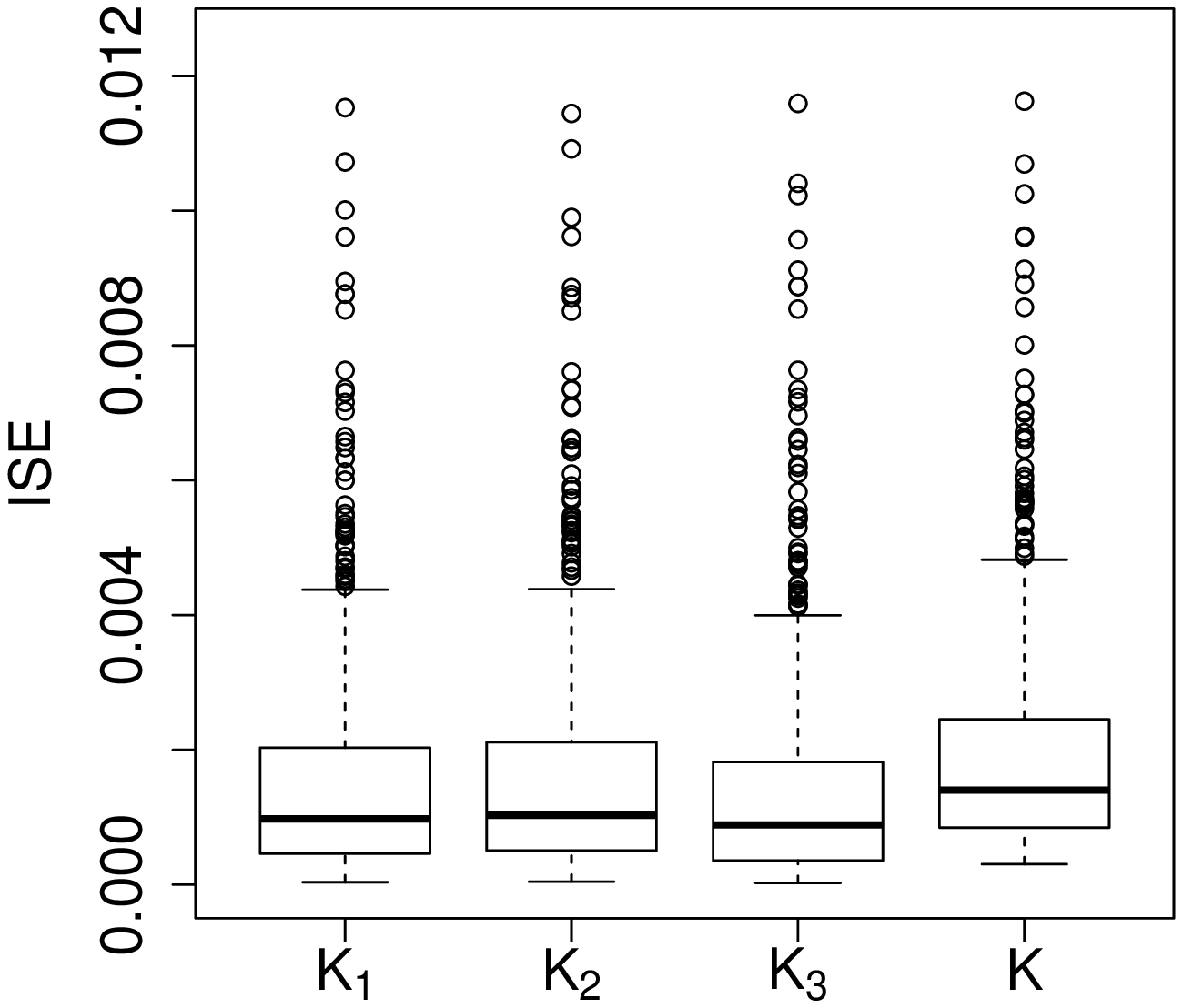}
\end{tabular}

\vspace{-1ex}

\caption{\normalsize\it $\mathrm{ISE}$ distributions for the boundary corrected
estimators with left boundary kernels $K^L_q$, $q=1,2,3$, and for the classical estimator
with kernel $K$ over the regions $[0,h]$ (left), $[0,1-h]$ (center) and $[0,1]$ (right).
$F$ is the beta mixture distribution $w B(1,2)+(1-w) B(2,b)$ with
$F^{\prime}_+(0)=1.5$ and $F^{\prime\prime}_+(0)=6$.
The boxplots are based on 500 generated samples of size $n=50$ and $K$ is the Epanechnikov kernel.} \label{iseboundary}
\end{figure}

We finish this section with a cautionary note that aims to call the attention of the reader
to the fact that, due to the continuity of $F$ on $\R$,
the boundary effects for kernel distribution function estimation
may not have the same impact in the global performance of the estimator as in
probability density or regression function estimation frameworks \citep[see][]{GasM:79}.
However, one may have cases where the local behaviour dominates the global behaviour
of the estimator which stresses the relevance in using boundary corrections for
the classical kernel distribution function estimator. We illustrate this fact by taking
the above considered beta mixture distribution
with $F^{\prime}_+(0)=1.5$ and $F^{\prime\prime}_+(0)=6$ (see Figure \ref{testdist}).
In Figure \ref{iseboundary} we present the empirical distribution of the integrated square error
of the classical estimator with kernel $K$ and of the boundary corrected estimators
with boundary kernels $K^L_q$, $q=1,2,3$,
over the boundary regions $[0,h]$ (left boundary ISE) and $[1-h,1]$ (right boundary ISE),
and over the all interval $[0,1]$ (ISE).
The boxplots are based on 500 generated samples of size $n=50$.
We conclude that the local behaviour of the estimator over the left boundary region
has a clear impact on the global performance of the estimator which supports the use of boundary corrections
for the classical kernel distribution function estimator.

\section{Proofs} \label{proofs}

We limit ourselves to present the proof of Theorem \ref{theoremboundL}.
The proofs of Theorems \ref{misetheo} and \ref{supremum} follow straightforward from
the proofs of the corresponding results given in \citet{Ten:13}
and the asymptotic expansions for bias and variance of $\tF_{nh}(x)$ we present below.

\medskip

\noindent
\textbf{Proof of Theorem \ref{theoremboundL}.a):}
For
$x\in \,]a,a+h[$, the expectation of $\tF_{nh}(x)$ is given by
\begin{align*}
\E \tF_{nh}(x)
& = \int F(x-uh) K^L(u;(x-a)/h) \,du,
\end{align*}
\citep[see][p.~186]{Ten:13}.
By the continuity of the second derivative of $F$
on $[a,b]$ and Taylor's formula, we have
\begin{align} \label{taylor}
F(x-uh) = F(x) - u h F^\prime(x) + u^2 h^2 \int_0^1 (1-t) F^{\prime\prime}(x-tuh) \, dt,
\end{align}
for $-1 \leq u \leq (x-a)/h$ ,
from which we deduce that
\begin{align} \label{eqbias}
\E \tF_{nh}(x) - F(x) - \frac{h^2}{2} F^{\prime\prime}(x) \mu_{L}((x-a)/h) = A(x,h) + B(x,h),
\end{align}
where
\begin{align*}
A(x,h) & = F(x) \big( \mu_{0,L}((x-a)/h) - 1 \big) - hF^\prime(x) \mu_{1,L}((x-a)/h)  \\
& \quad + \frac{h^2}{2} F^{\prime\prime}(x) ((x-a)/h) \mu_{1,L}((x-a)/h),
\end{align*}
and
$$B(x,h) = h^2 \iint_0^1 (1-t) \big( F^{\prime\prime}(x-tuh) - F^{\prime\prime}(x) \big) dt\,  u^2 K^L(u;(x-a)/h) \,du,$$
is such that
\begin{align} \label{supB}
\sup_{x\in\,]a,a+h[} |B(x,h)| \leq \frac{h^2}{2} \sup_{\alpha\in\, ]0,1[} |\mu_{0,L}|(\alpha)
\sup_{y,z \in [a,b]:\,|y-z| \leq h} |F^{\prime\prime}(y)-F^{\prime\prime}(z)|.
\end{align}

On the other hand, taking into account that $F(a)=0$ and using condition (\ref{C2}) and the Taylor's expansions
\begin{align} \label{taylor1}
F(x) & = (x-a) F^\prime(a) + \frac{1}{2} (x-a)^2 F^{\prime\prime}(a) \nonumber \\
& \quad + (x-a)^2 \int_0^1 (1-t) \big( F^{\prime\prime}(a+ (x-a)t ) - F^{\prime\prime}(a) \big) dt
\end{align}
and
\begin{align} \label{taylor2}
F^\prime(x) = F^\prime(a) + (x-a) F^{\prime\prime}(a) +
(x-a) \int_0^1 \big( F^{\prime\prime}(a+ (x-a)t) - F^{\prime\prime}(a) \big) dt,
\end{align}
we get
\begin{align} \label{supA}
\sup_{x\in\,]a,a+h[} |A(x,h)| \leq h^2 \sup_{\alpha\in\, ]0,1[} |\mu_{0,L}|(\alpha)
\sup_{y,z \in [a,b]:\,|y-z| \leq h} |F^{\prime\prime}(y)-F^{\prime\prime}(z)|.
\end{align}

Part a) of Theorem \ref{theoremboundL} follows now from (\ref{eqbias}), (\ref{supB}) and (\ref{supA}), and the fact that
\begin{equation*}
\sup_{y,z \in [a,b]:\,|y-z| \leq h} |F^{\prime\prime}(y)-F^{\prime\prime}(z)|=o(1). \halmoseq
\end{equation*}

\bigskip

\noindent
\textbf{Proof of Theorem \ref{theoremboundL}.b):}
From Part a), the variance of $\tF_{nh}(x)$ is given by
\begin{align}
n \Var \tF_{nh}(x)
& = \int \bK^L(z;(x-a)/h)^2 hf(x-uh) dz - \big( \E \tF_{nh}(x) \big)^2  \nonumber \\ 
& = F(x)(1-F(x)) + C(x,h) + O\big( h^2 \big), \nonumber
\end{align}
uniformly in $x\in \,]a,a+h[$, where
$$C(x,h) = \int \bK^L(u;(x-a)/h)^2 hf(x-uh) du - F(x).$$
Moreover, using (\ref{taylor}) and the fact that
$$\lim_{u \rightarrow -\infty} \bK^L(u;\alpha) = 0 \quad \mbox{and }
\lim_{u \rightarrow +\infty} \bK^L(u;\alpha) = \mu_{0,L}(\alpha), \mbox{ for } \alpha\in\,]0,1[,$$
we deduce that
\begin{align}
C(x,h)
& = \int F(x-zh)B^L(z;(x-a)/h) dz - F(x) \nonumber \\
& = F(x) \big( \mu_{0,L}((x-a)/h)^2 - 1 \big) - h F^\prime(x) m_{1,L}((x-a)/h) \nonumber \\
& \quad + h^2 \iint_0^1 (1-t) F^{\prime\prime}(x-tuh) dt u^2 B^L(u;(x-a)/h) du \nonumber \\
& = F(x) \big( \mu_{0,L}((x-a)/h)^2 - 1 \big) - h F^\prime(x) m_{1,L}((x-a)/h) + O(h^2), \label{expanC}
\end{align}
uniformly in $x\in \,]a,a+h[$, as $\sup_{\alpha\in\,]0,1[} \int |u^2 B^L(u;\alpha)| du < \infty$.

Finally, from (\ref{expanC}) and Taylor's expansions (\ref{taylor1}) and (\ref{taylor2})
we get
$$\sup_{x\in\,]a,a+h[} \left| C(x,h) + h \, F^\prime(x) \nu_L\big((x-a)/h\big) \right| = O(h^2),$$
which concludes the proof.
\hfill$\blacksquare$

\bigskip\bigskip

\noindent{\bf Acknowledgments.}
Research partially supported by
Centro de Matem\'atica da Universidade de
Coimbra (funded by the European Regional Development Fund through the program COMPETE and by the
Portuguese Government through the FCT--Funda\c{c}\~ao para a Ci\^encia e Tecnologia under the
project PEst-C/MAT/UI0324/2011).

\end{document}